\documentclass[fleqn]{mat01}
\usepackage{times,mathtimy,amssymb,latexsym,boxit}
\begin{document}

\setcounter{page}{251}
\firstpage{251}

\font\xx=msam5 at 10pt
\def\ab{\mbox{\xx{\char'03}}}

\newcommand{\Lam}{\Lambda}
\newcommand{\lam}{\lambda}
\newcommand{\tet}{\theta}
\newcommand{\sym}[3]{#1^{#2_{#3}}}
\newcommand{\til}[1]{\widetilde{#1}}
\newcommand{\innp}[2]{\left\langle #1,#2\right\rangle}
\newcommand{\summ}[3]{\sum_{#1=#2}^{#3}}
\newcommand{\cdel}{c~-~\frac{1}{2\pi}[(\rho_{g}^{-1}\delta)_{{\cal H}_{g}}]}
\newcommand{\inte}[2]{\int_{#1}^{#2}}
\newcommand{\kernel}[1]{\mbox{Ker}\,{#1}}
\newcommand{\coker}[1]{\mbox{Coker}\,{#1}}
\newcommand{\codim}{\mbox{codim}}
\newcommand{\im}[1]{\mbox{Im}\,{#1}}
\newcommand{\norml}{\left\|}
\newcommand{\normr}{\right\|}
\newcommand{\brce}[2]{_{\{#1,#2\}}}
\newcommand{\fr}[2]{\frac{#1}{#2}}
\newcommand{\alp}{\alpha}
\newcommand{\bet}{\beta}
\newcommand{\gam}{\gamma}
\newcommand{\del}{\delta}
\newcommand{\eps}{\epsilon}
\newcommand{\sig}{\sigma}
\newcommand{\smy}[3]{#1_{#2}^{#3}}
\newcommand{\transv}{\frown\!\!\!\!\mid\,}

\newcommand{\bol}[2]{{\bf #1}^{#2}}
\newcommand{\gras}[3]{G_{#1}({\Bbb #2}^{#3})}
\newcommand{\rarr}{\rightarrow}
\newcommand{\larr}{\leftarrow}
\newcommand{\mtxii}[4]{\left(\begin{array}{cc}#1&#2\\#3&#4\end{array}\right)}
\newcommand{\R}{{\Bbb R}}
\newcommand{\C}{{\Bbb C}}
\newcommand{\Q}{{\Bbb Q}}
\newcommand{\Z}{{\Bbb Z}}
\newcommand{\N}{{\Bbb N}}
\newcommand{\ov}{\overline}
\newcommand{\rot}[1]{#1^{\frac{1}{2}}}
\newcommand{\iroot}[1]{#1^{-\frac{1}{2}}}
\newcommand{\noin}{\noindent}
\newcommand{\zc}{\bar{z}}
\newcommand{\tr}{\mbox{tr}}
\newcommand{\fun}[2]{\pi_{1}(#1,#2_{0})}
\newcommand{\pair}[3]{(#1, #2_{#3})}
\newcommand{\boun}{\partial}
\newcommand{\ext}{\bigwedge}
\newcommand{\affine}[2]{{\Bbb A}^{#1}(#2)}
\newcommand{\proj}[2]{{\Bbb P}^{#1}(#2)}
\newcommand{\itm}[1]{\item[{\rm (#1)}]}
\newcommand{\eucl}[2]{{\Bbb {#1}}^{#2}}
\newcommand{\project}[2]{{\Bbb #1}{\Bbb P}(#2)}
\newcommand{\ad}[1]{{\mbox{ad}{\,#1}}}
\newcommand{\Ad}[1]{{\mbox{Ad}{\,#1}}}
\newcommand{\grad}[1]{{\mbox{grad}{\,#1}}}
\newcommand{\Exp}{\mbox{Exp}}

\newtheorem{theo}{Theorem}
\renewcommand\thetheo{\arabic{section}.\arabic{theo}}
\newtheorem{theor}[theo]{\bf Theorem}
\newtheorem{lem}[theo]{Lemma}
\newtheorem{propo}[theo]{\rm PROPOSITION}
\newtheorem{rmk}[theo]{Remark}
\newtheorem{defn}[theo]{\rm DEFINITION}
\newtheorem{exam}[theo]{Example}
\newtheorem{coro}[theo]{\rm COROLLARY}
\newtheorem{nota}[theo]{\it Notation}
\def\conjecture{\trivlist\item[\hskip\labelsep{\it Conjecture.}]}
\def\claim{\trivlist\item[\hskip\labelsep{\it Claim.}]}
\def\poclaim{\trivlist\item[\hskip\labelsep{\it Proof of Claim.}]}

\title{Subanalytic bundles and tubular neighbourhoods of zero-loci}

\markboth{Vishwambhar Pati}{Subanalytic bundles and tubular neighbourhoods of zero-loci}

\author{VISHWAMBHAR PATI}

\address{Stat-Math Unit, Indian Statistical Institute, R.V. College Post,
Bangalore~560~059, India\\
\noindent E-mail: pati@isibang.ac.in}

\volume{113}

\mon{August}

\parts{3}

\Date{MS received 18 December 2002}

\begin{abstract}
We introduce the natural and fairly general notion of a subanalytic
bundle (with a finite dimensional vector space $P$ of sections) on a
subanalytic subset $X$ of a real analytic manifold $M$, and prove that
when $M$ is compact, there is a Baire subset $U$ of sections in $P$
whose zero-loci in $X$ have tubular neighbourhoods, homeomorphic to the
restriction of the given bundle to these zero-loci.
\end{abstract}


\maketitle

\section{Introduction}

In this paper, we introduce the notion of a subanalytic bundle $E$
(generated by a finite dimensional space $P$ of global sections) on a
(not necessarily closed) subanalytic set $X$ inside a real analytic
manifold $M$, as a natural generalisation of real analytic bundles on
real analytic spaces to the subanalytic setting. We prove (in
Theorem~\ref{strattubular2} below) that for $M$ compact, there exists a
Baire subset $U$ of sections in $P$, such that for $s\in U$, there exist
tubular neighbourhoods of the zero-locus $Z=s^{-1}(0_{E})$ of $s$ in
$X$, i.e. which are homeomorphic to the restriction of the given bundle
to $Z$. To keep the account self-contained we recall basic facts about
subanalytic sets in \S2 and Strong Whitney (SW) stratifications
(defined by Verdier) in \S4.

\looseness 1 We remark here that the main Theorem \ref{strattubular2} would follow
from Theorem~1.11 on p.~48 of [G-M]. However, the proof (`deformation
to the normal bundle') sketched in [G-M] is incomplete, at least in the
generality that it is stated. In this generality, the stratified
submersion they construct is not proper (as was pointed out by
V~Srinivas), and hence Thom's First Isotopy Lemma is inapplicable. To
circumvent this, we have imposed the hypothesis of compactness on the
ambient real analytic manifold $M$ containing the subanalytic set $X$,
but no compactness assumption on $X$. Our hypotheses are general enough
to cover most situations arising in real or complex algebraic geometry
(see Example~\ref{subanalyticexample} and Remark~\ref{strattubular3}).

\section{Subanalytic sets and maps}

Let $M$ be a real-analytic manifold. We will always assume $M$ to be
connected, Hausdorff, second countable and paracompact.

\begin{defn}$\left.\right.$\vspace{.5pc}

\noindent {\rm We say $X\subset M$ is a {\em subanalytic set} of $M$
if there exists an open covering
${\cal U}$ of $M$ ({\em not} just of $X$) such that for each $U\in {\cal U}$,
\begin{equation*}
X\cap U = \bigcup_{i=1}^{p}\left(f_{i1}(A_{i1})-f_{i2}(A_{i2})\right),
\end{equation*}
where $f_{ij}:N_{ij}\rarr U$ for $1\leq i\leq p$ and $j=1,2$, are real analytic maps of real analytic
manifolds $N_{ij}$, $A_{ij}$ are closed analytic subsets of $N_{ij}$ and $f_{ij|A_{ij}}$
are proper maps (see Proposition~3.13 in [B-M] and Definition~3.1 in [Hi]).\label{subanaldef}}
\end{defn}

\begin{exam}{\rm All real (resp. complex) analytic subsets of a real
(resp. complex) analytic manifold are subanalytic sets. In particular, (real or
complex) algebraic subsets of a (real or complex) algebraic manifold (such as
projective space, or Grassmannians) are subanalytic sets. Also since
subanalytic subsets of a real analytic manifold form a Boolean algebra (see (i) of
Proposition~\ref{subanalfacts} below) all real (resp. complex) analytically (or algebraically) constructible
sets in a real (resp. complex) analytic (or algebraic) manifold are
subanalytic sets. In particular, all (real or complex) affine algebraic
varieties are subanalytic in both affine space, and projective space. Real or
complex quasiprojective varieties are subanalytic sets in the
corresponding projective spaces.\label{subanalyticexample}}
\end{exam}

\begin{rmk}
{\rm A real analytic subset $X$ (subspace) of a real analytic manifold $M$ is a
closed subset of $M$ by definition. In particular if $M$ is compact, so is $X$.
By contrast, a subanalytic set $X$ of a real analytic manifold $M$ need not
be closed, and need not be compact even if $M$ is compact.
\label{analsubanal}}
\end{rmk}

\begin{defn}$\left.\right.$\vspace{.5pc}

\noindent {\rm Let $X\subset M$
and $Y\subset N$ be subanalytic sets in the real analytic manifolds
$M, N$ respectively. We say that a map $f:(X,M)\rarr (Y,N)$
is a {\em subanalytic map} if $f:X\rarr Y$ is a continuous map, and the graph
\begin{equation*}
\Gamma_{f}:=\{(x,y)\in M\times N:\, x\in X, \,\,y=f(x)\}
\end{equation*}
is a subanalytic set in $M\times N$ (see [Ha], 4.1, or Definition~3.2 in [B-M]). Note that although
the map $f$ is defined only on $X$, its subanalyticity depends on the
ambient $M$, $N$, as we shall see in Remark~\ref{subanalrmk} below.
\label{subanalmapdef}}
\end{defn}

\begin{nota}{\rm If $X,M,N$ are as above, and $f:(X,M)\rarr (N,N)$ is a
subanalytic map, we shall write $f:(X,M)\rarr N$ is a subanalytic map,
for notational convenience.}
\end{nota}

\begin{rmk}
{\rm The subanalyticity (or analyticity) of a set, or of a map {\em depends on the
ambient spaces} $M, N$. For example, $X=\{\fr{1}{n}\}_{n\in {\Bbb N}}$ is a subanalytic
set in  $(0,\infty)$, but not in $\R$. In the former, it is the zero set
of the analytic function $\sin\fr{\pi}{x}$, so analytic and hence
subanalytic in $(0,\infty)$. It is not subanalytic in $\R$ because the connected components
of its germ at $0$ in $\R$ do not form a locally finite collection (see (viii) of Proposition
\ref{subanalfacts} below).

Similarly, the map $((0,1), (0,1))\rarr (\R, \R)$ defined by $x\mapsto
\sin\,\fr{\pi}{x}$ is clearly subanalytic, because its graph
$\Gamma:=\{(x,\sin\,\fr{\pi}{x}):x\in (0,1)\}$ is an analytic (hence
subanalytic) subset of $(0,1)\times\R$. On the other hand, the same
mapping regarded as a map $((0,1), \R)\rarr (\R, \R)$ is {\em not}
subanalytic, since $\Gamma$ is not a subanalytic subset in $\R\times
\R$. (By (i) of Proposition~\ref{subanalfacts}, if $\Gamma$ were
subanalytic in $\R\times \R$, its intersection with the $x$-axis would
have to be subanalytic in $\R\times\R$. But this intersection is the set
$\{(\fr{1}{n},0)\}_{n\in {\Bbb N}}$, which is not subanalytic because the connected components of its
germ at $(0,0)$ in $\R\times\R$ is not a locally finite collection (see (viii) of
Proposition~\ref{subanalfacts} below.).\label{subanalrmk}}
\end{rmk}

\begin{propo}{\rm (Facts on subanalytic sets and maps)}$\left.\right.$\vspace{.5pc}

\noindent We collect some well-known facts on subanalytic sets and maps{\rm :}

\begin{enumerate}
\renewcommand{\labelenumi}{\rm (\roman{enumi})}
\leftskip 1pc
\item The collection of subanalytic sets of a real
analytic manifold $M$ forms a Boolean algebra.

\item If $f:M\rarr N$ is a proper real analytic map of real analytic
manifolds{\rm ,} and $X\subset M$ a subanalytic set{\rm ,} then $f(X)$ is
subanalytic in $N$. In particular{\rm ,} if $M$ is compact{\rm ,}
the image $f(X)$ is subanalytic. If $g:(X,M)\rarr N$ is a subanalytic map
and $X$ is relatively compact in $M${\rm ,} then $g(X)\subset N$ is
subanalytic.

\item If $X\subset M$ and $Y\subset N$ are subanalytic{\rm ,} then
$X\times Y$ is subanalytic in $M\times N$.

\item If $X\subset M$ is subanalytic{\rm ,} then the diagonal
\begin{equation*}
\Delta_{X}:=\{(x,x)\in M\times M: \;x\in X\}
\end{equation*}
is subanalytic in $M\times M$. Thus the inclusion map $i:(X,M)\rarr M$
is always a subanalytic map.

\item Let $f:M\rarr N$ be a real analytic map of real analytic
manifolds. If $X\subset M$ is a subanalytic set{\rm ,} then the restricted
map $f_{|X}:(X,M)\rarr (N,N)$ is a subanalytic map.

\item Let $f:M\rarr N$ be a real analytic map of real analytic
manifolds. If $Y\subset N$ is subanalytic{\rm ,} then $f^{-1}(Y)\subset M$ is
subanalytic.

\item The closure $\ov{X}$ of a subanalytic set $X\subset M$ is also
subanalytic.

\item Let $X\subset M$ be a subanalytic set in a real analytic
manifold $M$. Then each connected component of $X$ is also subanalytic.
The collection of connected components of $X$ is a locally finite
collection in $M$.
\end{enumerate}\label{subanalfacts}
\end{propo}

\begin{proof}
For a proof of (i), see Proposition~3.2 of [Hi] or \S3 of [B-M].

For a proof of the first statement of (ii), see  [Hi], Proposition~3.8. For a
proof of the second statement of (ii), see the remark after Definition~3.2 in
[B-M]. Easy examples can be constructed to show that the properness condition
cannot be dropped from the\break hypothesis.

To see (iii), first note that $X\subset M$ subanalytic implies $X\times
N\subset M\times N$ is also subanalytic. For, in the Definition
\ref{subanaldef} above, one merely takes the open covering
${\cal U}\times N:=\{U\times N:\;U\in {\cal U}\}$, the closed
sets $A_{ij}\times N\subset N_{ij}\times N$ and the maps
$f_{ij}\times \mbox{id}_{N}:N_{ij}\times N\rarr U\times N$,
which are proper on $A_{ij}\times N$ since $f_{ij}$ are proper on $A_{ij}$.
Similarly $M\times Y\subset
M\times N$ is also subanalytic. By (i) above, the intersection
$X\times Y=(X\times N)\cap (M\times Y)$ is also subanalytic. This proves
(iii).

To see (iv), note that by (iii) $X\times X\subset M\times M$ is
subanalytic. The diagonal $\Delta_{M}\subset M\times M$ is subanalytic
since it is analytic in $M\times M$. The intersection
$\Delta_{X}=\Delta_{M}\cap (X\times X)$ is therefore subanalytic by (i).
Since $\Delta_{X}$ is the graph of the inclusion $i:(X,M)\rarr M$ inside $M\times M$, it
follows that $i$ is a subanalytic map.

To see (v), we note that the graph of $f_{|X}$ in $M\times N$ is just
the intersection of the graph $\Gamma_{f}$ of $f$ and $X\times N$ inside $M\times N$. Since
$\Gamma_{f}$ is an analytic set in $M\times N$, it is subanalytic, and
since by (iii) $X\times N$ is subanalytic, their intersection is
subanalytic by (i).

For (vi), let $Y\subset N$ be subanalytic and
let ${\cal U}$ be an open covering of $N$ such that for each $U\in {\cal
U}$ we have
\begin{equation*}
Y\cap U= \bigcup_{i=1}^{p}\left(f_{i1}(A_{i1})-f_{i2}(A_{i2})\right),
\end{equation*}
where $f_{ij}:N_{ij}\rarr U$ are real analytic maps of real analytic
manifolds $N_{ij}$, $A_{ij}$ are closed analytic subsets of $N_{ij}$ and
$f_{ij|A_{ij}}$ are proper maps. Now take the open covering
\begin{equation*}
f^{-1}({\cal U}):=\{f^{-1}(U):U\in {\cal U}\}
\end{equation*}
of $M$, and set $\til{N}_{ij}:=N_{ij}\times f^{-1}(U)$, with
$\til{f}_{ij}:\til{N}_{ij}\rarr f^{-1}(U)$ being the second projection.
Let $\til{A}_{ij}=A_{ij}\times_{U}\,f^{-1}(U)$ (the fibre product, a
closed analytic subset of $\til{N}_{ij}$). Observe that the restriction
to $\til{A}_{ij}$ of the natural real analytic projection
$\til{f}_{ij}:\til{N}_{ij}\rarr f^{-1}(U)$, is the `base change' to
$f^{-1}(U)$ of the restriction $f_{ij|A_{ij}}:A_{ij}\rarr U$, and this
last map is given to be proper. Hence this restriction
$\til{f}_{ij|\til{A}_{ij}}$ is proper. It is easily verified that
\begin{equation*}
f^{-1}(Y)\cap f^{-1}(U)=
\bigcup_{i=1}^{p}\left(\til{f}_{i1}(\til{A}_{i1})-\til{f}_{i2}(\til{A}_{i2})
\right)
\end{equation*}
which shows that $f^{-1}(Y)$ is subanalytic.

For (vii), see the immediate consequences following Definition~3.1 in [B-M], and
also Corollary~3.2.9 in [Hi].

For (viii), see the immediate consequences following Definition~3.1 in [B-M].
(Also see Proposition~3.6 and Corollary~3.7.10 in [Hi].)\hfill \ab
\end{proof}

For an analytic subset $X$ inside a real analytic manifold $M$,
there is a structure sheaf, making it a locally
ringed space. Thus mappings ($=$ morphisms) of real analytic spaces are
easy to define, and obey the usual functorial properties. For
subanalytic sets inside a real analytic manifold, we note that there is
no such structure sheaf, and the definition of a subanalytic map is
dependent on the ambient manifold $M$. Thus the notion of
`subanalytic equivalence' of subanalytic sets $X\subset M$ and $Y\subset N$
requires some care. We propose one such below, which may not be the most
general, but is good enough for our purposes. We are unaware if this notion
exists in the literature.

\begin{lem} Let $M${\rm ,} $M_{1}$ be real analytic manifolds. Let $X\subset M$ be a subanalytic set{\rm ,} and suppose
$j:(X,M)\rarr M_{1}$ is a subanalytic map. Suppose there exists a
{\em proper} real analytic map $p:M_{1}\rarr M$ such that
$p\circ j=\mbox{id}_{X}$. Then

\begin{enumerate}
\renewcommand{\labelenumi}{\rm (\roman{enumi})}
\leftskip .5pc
\item $j:X\rarr j(X)$ is a homeomorphism{\rm ,} and its image $X_{1}:=j(X)\subset M_{1}$
is subanalytic in $M_{1}$. Further $X\subset M$ is relatively compact if and only if
$X_{1}=j(X)\subset M_{1}$ is relatively compact.

\item For each subanalytic map $f:(X_{1},M_{1})\rarr N, N$
a real analytic manifold{\rm ,} the composite map{\rm :}
\begin{equation*}
\hskip -1.3pc f\circ j:(X,M)\rarr N
\end{equation*}
is a subanalytic map.

\item For each subanalytic map $g:(X,M)\rarr N, N$ a real analytic manifold{\rm ,}
the composite:
\begin{equation*}
\hskip -1.3pc g\circ p_{|X_{1}}:(X_{1},M_{1})\rarr N
\end{equation*}
is a subanalytic map.
\end{enumerate}\label{properlemma}
\end{lem}

\begin{proof}
It is clear that $j$ is a homeomorphism, with inverse $p_{|j(X)}$.
Consider the real analytic map
\begin{align*}
\theta: M_{1} &\rarr M\times M_{1}\\[.2pc]
z &\mapsto (p(z),z).
\end{align*}
Also let
\begin{equation*}
\Gamma_{j}=\{(x,j(x))\in M\times M_{1}:\,x\in X\}
\end{equation*}
be the graph of $j$. Since $j$ is a subanalytic map, $\Gamma_{j}\subset
M\times M_{1}$ is subanalytic.

We claim that the inverse image $\theta^{-1}(\Gamma_{j})\subset M_{1}$ is
precisely $X_{1}$. For $z\in \theta^{-1}(\Gamma_{j})\;\Rightarrow
\;(p(z),z)\in\Gamma_{j}\;\Rightarrow\;p(z)\in X$ and
$z=jp(z)\;\Rightarrow\;z\in j(X)$. Conversely,  $z=j(x)$ for $x\in
X\;\Rightarrow\;\theta(z)=(p\circ j(x), j(x))=(x,j(x))$ which is clearly
in $\Gamma_{j}$. Hence the claim. Now, since $\theta$ is real analytic, and $\Gamma_{j}$
is subanalytic, we have by (vi)
of Proposition \ref{subanalfacts} that $X_{1}=\theta^{-
1}(\Gamma_{j})\subset M_{1}$ is subanalytic. This proves the first assertion of
(i). For the second assertion, note that the continuity of $p$ implies
$p(\ov{X}_{1})\subset \ov{p(X_{1})}=\ov{X}$. Since $p$ is proper, it
is a closed map, and so $p(\ov{X}_{1})$ is a closed set containing $p(X_{1})=X$.
Hence $\ov{X}=p(\ov{X}_{1})$. Thus if $X_{1}$ is relatively compact in $M_{1}$,
$X$ is relatively compact in $M$. Conversely if $X$ is relatively compact
in $M$, $\ov{X}_{1}\subset p^{-1}(\ov{X})$ and $p$ is proper implies
that $\ov{X}_{1}$ is a closed subset of the compact set $p^{-
1}(\ov{X})$, and hence also compact. That is, $X_{1}$ is relatively compact in
$M_{1}$. This proves (i).

To see (ii), let $f:(X_{1},M_{1})\rarr N$ be a subanalytic map. Thus the
graph
\begin{equation*}
\Gamma_{f} = \{(x,f(x))\in M_{1}\times N:\; x\in X_{1}\}
\end{equation*}
is a subanalytic set. The real analytic map:
\begin{equation*}
(p\times \mbox{id}):M_{1}\times N \rarr M\times N
\end{equation*}
is proper, since $p$ is proper. But the image
\begin{align*}
(p\times \mbox{id})\left(\Gamma_{f}\right) &= \{(p(x),f(x))\in M\times N:\;x\in
X_{1}\}\\[.2pc]
&= \{(pj(y),fj(y))\in M\times N:\;y\in X\}\\[.2pc]
&= \{(y,fj(y))\in M\times N:\;y\in X\}\\[.2pc]
&= \Gamma_{f\circ j},
\end{align*}
the graph of the composite $f\circ j$. Since $p\times \mbox{id}$ is
a proper real analytic map, and $\Gamma_{f}$ is subanalytic, it follows by (ii) of
Proposition~\ref{subanalfacts} that this image, the graph
$\Gamma_{f\circ j}$ is subanalytic in $M\times N$, so that $f\circ
j:(X,M)\rarr N$ is a subanalytic map. This\break proves (ii).

To see (iii), let $g:(X,M)\rarr N$ be a subanalytic map. This means that
the graph $\Gamma_{g}\in M\times N$ is a subanalytic set. Consider the set
\begin{align*}
\hskip -1pc (p\times \mbox{id})^{-1}(\Gamma_{g})\cap (X_{1}\times N) &= \{(m,n)\in M_{1}\times N:
\;(p(m),n)\in \Gamma_{g}, \; m\in X_{1}\}\\[.2pc]
\hskip -1pc &= \{(m,n)\in M_{1}\times N:\; n=gp(m), \;m\in X_{1}\}\\[.2pc]
\hskip -1pc &= \Gamma_{g\circ p_{|X_{1}}}.
\end{align*}
Since $X_{1}\times N$ is subanalytic in $M_{1}\times N$ by (iii) of
Proposition \ref{subanalfacts},
$\Gamma_{g}$ is subanalytic in $M\times N$ by definition,
and $p\times\mbox{id}:M_{1}\times N\rarr M\times N$ is real
analytic, (iii) follows from (i) and (vi) of Proposition
\ref{subanalfacts}.\hfill \ab
\end{proof}

The above Lemma~\ref{properlemma} shows that under the hypotheses stated
there, the subanalytic sets $(X,M)$ and $(X_{1},M_{1})$ are `equivalent'
in some sense. More precisely, we make the following definition:

\begin{defn}{\rm (Pseudoequivalence of subanalytic sets)$\left.\right.$\vspace{.5pc}

\noindent Let $M,\;M_{1}$ be real
analytic manifolds, with $X\subset M$ a subanalytic set
and $j:(X,M)\rarr M_{1}$ a subanalytic map. If there
exists a {\em proper real analytic map} $p:M_{1}\rarr M$ such that
$p\circ j=\mbox{id}_{X}$, then we say that the subanalytic sets $(X,M)$
and $(j(X),M_{1})$ are {\em subanalytically pseudoequivalent}. The map $j$ is called
a {\em subanalytic pseudoequivalence}. We note that $X$ and $X_{1}:=j(X)$ are
therefore {\em a fortiori} homeomorphic, and also (i) of the
Lemma~\ref{properlemma} implies that $X$ is relatively compact in $M$ iff $X_{1}$ is relatively
compact in $M_{1}$. \label{properequivdef}}
\end{defn}

The prototypical example of such a subanalytic pseudoequivalence of interest to us
in the sequel is the following.

\begin{exam}{\rm (Graph embeddings).\ \ \ Let $X\subset M$ be a subanalytic set,
and $f:(X,M)\rarr N$ a subanalytic map. {\em Assume that $N$ is
compact}. Set $M_{1}:=M\times N$, and let $j:(X,M)\rarr M_{1}$ be the
graph embedding defined by $j(x)=(x,f(x))$ for $x\in X$. $j$ is a
subanalytic map because its graph in $M\times M_{1}$ is the set
$\{(x,x,f(x)):\,x\in X\}$, which is
precisely the intersection of the two subanalytic sets
$\Delta_{X}\times N$ and $M\times \Gamma_{f}$ in $M\times M_{1}$,
and therefore subanalytic (by Proposition~\ref{subanalfacts}(i), (iii), and (iv)).
The projection $p:M\times N\rarr M$ is proper
since $N$ is compact, and thus we have the requirements of Definition
\ref{properequivdef}. That is, $j:(X,M)\rarr (\Gamma_{f}, M\times N)$ is
a subanalytic pseudoequivalence.\label{graphexample}}
\end{exam}

\section{Analytic bundle theory}

\setcounter{theo}{0}

We review some basic notions of bundles from the real-analytic set-up,
with a view to generalising them to the subanalytic set-up.

Suppose that $X\subset M$ is a real analytic subset ($=$ subspace) in a
real analytic manifold $M$. By definition, its germ at each point of $M$
(not just $X$) is given by the vanishing of some ideal, so by definition
$X$ is closed in $M$. Then $X$ comes equipped with a structure sheaf
${\cal O}_{X}$, whose stalk ${\cal O}_{X,x}$ at $x\in X$ consists of
germs of real analytic functions on $X$ at $x$. It is, by definition,
the local ring ${\cal O}_{M,x}/{\cal I}_{X,x}$, where ${\cal O}_{M,x}$
is the local ring of $M$ at $x$ consisting of convergent power series at
$x\in M$, and ${\cal I}_{X,x}$ is the ideal of functions in ${\cal
O}_{M,x}$ vanishing on the germ of $X$ at $x$.

Let $\pi:E\rarr X$ be a {\em real analytic vector bundle of rank} $k$ on
$X$. That is, its transition cocycles $g_{ij}:U_{i}\cap U_{j}\rarr
GL(k,\R)$ are real-analytic functions for all $i,j$. The sheaf (of germs
of analytic sections) of an analytic vector bundle $E$ on $X$ is a
locally free sheaf ${\cal E}$ of modules over the structure sheaf ${\cal
O}_{X}$. Global sections of this sheaf are called global sections of
$E$. $E$ is said to be {\em generated by a vector space} $P\subset {\cal
E}(X)$ {\em of global sections} if the natural sheaf map:
\begin{equation*}
P\otimes {\cal O}_{X}\rightarrow {\cal E}
\end{equation*}
is a surjective map of ${\cal O}_{X}$-modules. This is equivalent to demanding that the
evaluation map:
\begin{align*}
\eps_{x}: P &\rightarrow E_{x}\\
s &\mapsto s(x)
\end{align*}
is surjective for all $x\in X$, where $E_{x}$ is the fibre of $E$ at $x$.

It is clear from the cocycle formulation of analytic vector bundles
that the natural bundle operations, such as direct sums, tensor
products, homs, duals and pullbacks under analytic maps of analytic
vector bundles are again analytic. Real analytic sections of the real
analytic bundle $\hom(E,F)$ are defined to be real analytic bundle
morphisms.

\begin{defn}{\rm (Universal exact sequence)$\left.\right.$\vspace{.5pc}

\noindent Let $G(n-k,P)$ denote the Grassmannian of
$(n-k)$-dimensional subspaces of $P$, where $P$ is a finite dimensional
real vector space of dimension $n$. On $G(n-k,P)$, there is the short
exact sequence of real analytic vector bundles and real analytic
morphisms:
\begin{equation}
0\rarr \gamma^{n-k}\rarr G(n-k,P)\times P\stackrel{\phi}{\rarr}
\nu^{k}\rarr 0,\label{univseq}
\end{equation}
where $G(n-k,P)\times P$ is the trivial rank-$n$ bundle on $G(n-k,P)$,
$\gamma^{n-k}$ is the tautological rank-$(n-k)$ real analytic bundle on
$G(n-k,P)$ (having fibre $V$ over the point $V\in G(n-k,P)$). The bundle
$\nu^{k}$ is the universal quotient bundle of rank $k$ on $G(n-k,P)$.
$P$ gets identified with the constant sections of $G(n-k,P)\times P$, and the
bundle $\nu^{k}$ is generated by the global sections $(\phi\circ s)$,
$s\in P$. \label{univdef}}
\end{defn}

Indeed, the bundles and morphisms defined above are all real algebraic, and
hence real analytic.

\begin{lem}{\rm (}Classifying maps{\rm )}.\ \ \ Let $M$ be a real analytic
manifold{\rm ,} and $X\subset M$ be a real analytic subset {\rm
(}$=$ subspace{\rm )}. Let $\pi:E\rarr X$ be a real analytic vector
bundle of rank $k$ with corresponding sheaf ${\cal E}$. Let $P$ be an
$n$-dimensional real vector subspace of ${\cal E}(X)$. Then the
following are equivalent{\rm :}
\begin{enumerate}
\renewcommand{\labelenumi}{\rm (\roman{enumi})}
\leftskip .5pc
\item $E\rarr X$ is generated by the global sections $P$.

\item There exists a real analytic map $f:X\rarr G(n-k, P)$ called the
{\em classifying map} such that the pullback of the universal short exact sequence $(1)$
under $f$ yields the short exact sequence
\begin{equation}
\hskip -1.3pc 0\rarr f^{\ast}\gamma^{n-k}\rarr X\times P\stackrel{\eps}{\rarr} E\rarr 0
\label{exactseq}
\end{equation}
on $X$. Here $\eps$ is the evaluation map $(x,s)\mapsto s(x)$.

\item There is a real analytic manifold $M_{1}${\rm ,} a real analytic
map $j:X\rarr M_{1}${\rm ,} and a proper real analytic map $p:M_{1}\rarr M$
such that{\rm :}

\leftskip .2pc
{\rm (a)} $p\circ j= {\rm id}_{X}$ and $X_{1}:=j(X)$ is isomorphic
to $X$ as a real analytic space via $j$.

{\rm (b)} There is an exact sequence of real analytic bundles on
$M_{1}${\rm :}
\begin{equation*}
\hskip -2.5pc 0\rarr K\rarr M_{1}\times P\stackrel{\eps_{1}}{\rarr} C\rarr 0
\end{equation*}
with $C$ generated by a space of global sections equal to
$P${\rm ,} and such that the pullback of the last two terms to $X$ via the
analytic isomorphism $j$ is
the morphism of bundles $X\times P \stackrel{\eps}{\rarr} E\rarr 0${\rm ,} defined by
evaluation {\rm (}i.e. $\eps(x,s)=\eps_{x}(s)\break =s(x)${\rm )}.

{\rm (}In keeping with Definition~\textup{\ref{properequivdef}}{\rm ,} one may want to call $j$
an {\em analytic pseudoequivalence}.{\rm )}

\end{enumerate}\label{analclass}
\end{lem}\vspace{-2pc}

\begin{proof}$\left.\right.$

\noindent (i) $\Rightarrow$ (ii). First we note that by
choosing a basis of $P$ of real analytic sections $\{e_{i}\}_{i=1}^{n}$,
that the maps:
\begin{align*}
\eps_{i}: X &\rarr E\\
x &\mapsto e_{i}(x)
\end{align*}
are all real analytic, so that the map:
\begin{align*}
\eps: X\times P &\rarr E\\
\left(x, s = \sum_{i}a_{i}e_{i}\right) &\mapsto \sum_{i}a_{i}\eps_{i}(x)
= s(x)
\end{align*}
is also real analytic. In particular, the map $\eps_{x}:P\rarr E$
defined by $\eps_{x}(s)=\eps(x,s)$ is also real analytic.

The classifying map $f$ is now defined as the map $x\mapsto
\ker\,\eps_{x}$, where $\eps_{x}:P\rarr E_{x}$ is the evaluation map. To
show that this map $f$ is real analytic, it is enough to do it locally.
Since $E$ is analytically locally trivial, we may assume without loss of
generality that $E$ is trivial. In this case, the map $f$ is just the
composite:
\begin{equation*}
X\rarr S\hom_{\R}(P,\R^{k})\rarr G(n-k, P),
\end{equation*}
where $S\hom_{\R}(P,\R^{k})$ is the open subset of $\hom_{\R}(P,\R^{k})$
consisting of surjective maps, the first arrow is $x\mapsto\eps_{x}$,
and the second arrow is the map $L\mapsto\ker L$. Since both these maps
are real analytic, the composite is real analytic. Thus the classifying
map $f:X\rarr G(n-k,P)$ is real analytic. So the bundle
$f^{\ast}(\gamma^{n-k})$ is a real analytic bundle, and it is clear that
its fibre is
\begin{equation*}
\left(f^{\ast}\gamma^{n-k}\right)_{x}=\ker{\left(\eps_{x}:P\rarr
E_{x}\right)}.
\end{equation*}
Finally, since $f^{\ast}(G(n-k,P)\times P)=X\times P$, and
$f^{\ast}\phi=\eps$, we have $f^{\ast}\nu^{k}=E$, and (ii) follows.

\noindent (ii) $\Rightarrow$ (iii). This follows by using the graph of the classifying
map $f:X\rarr G(n-k,P)$. More precisely, let us define
$M_{1}=M\times G(n-k,P)$, a real analytic manifold, and set
\begin{equation*}
X_{1}:=\{(x,y)\in M_{1}:\; y=f(x)\},
\end{equation*}
where $f$ is the classifying map of (ii). Define $j:X\rarr X_{1}$ to be
the graph embedding $j(x)=(x,f(x))$. Let $p:M_{1}\rarr M$ be the first
projection. Clearly $p\circ j=\mbox{id}_{X}$.

That $X_{1}$ is a real analytic space in $M_{1}$ is clear from the corresponding local fact,
i.e. that the germ of the graph $X_{1}$, at any point $(x,y)\in M_{1}$
is defined by the real analytic ideal ${\cal I}_{X_{1},(x,y)}$ generated by
${\cal I}_{X,x}\otimes 1$
and (the component functions of) $(1\otimes v)-(f(u)\otimes 1)$ in the
completed tensor product
\begin{equation*}
{\cal O}_{M_{1},(x,y)} = {\cal O}_{M,x}\;\widehat{\otimes}\;{\cal O}_{G(n-k,P),y},
\end{equation*}
where $u$ and $v$ are local coordinates around $x\in M$ and $y\in
G(n-k,P)$. The graph embedding:
\begin{align*}
j:X &\rarr M_{1}\\
x &\mapsto (x,f(x))
\end{align*}
is a real analytic map, with image $X_{1}$. Since the first projection
$p:M_{1}\rarr M$ provides the analytic inverse to $j:X\rarr X_{1}$, the
map $j$ is an analytic isomorphism of the real analytic spaces $X$ and
$X_{1}$. This proves (a).

If we let $p_{2}:M_{1}\rarr G(n-k,P)$ denote the real analytic map
defined by the second projection, and set $K:=p_{2}^{\ast}\gamma^{n-k}$,
and $C:=p_{2}^{\ast}\nu^{k}$, we have the $p_{2}$-pullback of the
universal exact sequence of (\ref{univseq}):
\begin{equation*}
0\rarr K\rarr M_{1}\times P\stackrel{\eps_{1}}{\rarr} C\rarr 0,
\end{equation*}
where $\eps_{1}:=p_{2}^{\ast}\phi$.
This is a short exact sequence of real analytic vector bundles on
$M_{1}$. If we pullback this sequence via the real analytic isomorphism
$j$, the fact that $p_{2}\circ j=f$ implies the short exact sequence
(\ref{exactseq}) on $X$. Thus we have the assertion (b).

\noindent (iii) $\Rightarrow$ (i). We pull back the  given exact sequence of real
analytic bundles on $M_{1}$ of (iii) to
$X$ via $j$. It continues to be an exact sequence of real analytic bundles
on $X$.
If we let $\eps$ denote the pullback morphism $j^{\ast}\eps_{1}$, and set
$j^{\ast}C=E$, the last two
terms of this pulled back exact sequence read:
\begin{equation*}
X\times P\stackrel{\eps}{\rarr}E\rarr 0
\end{equation*}
which shows that the $\eps$ images of the constant sections $x\mapsto
(x,s)$ generate $E$, and hence (i) follows.\hfill \ab
\end{proof}

\begin{rmk}
{\rm The equivalence of (i) and (iii) of the above lemma shows that in
considering analytic bundles generated by global sections on an analytic
subset $X\subset M$, we lose no generality in assuming (up to analytic
pseudoequivalence) that such bundles are {\em restrictions} of similar
(viz. generated by global sections) real analytic bundles on the {\em
ambient smooth} $M$ to $X$. We will generalise all this to a subanalytic
setting in the next section.\label{restrictrmk}}
\end{rmk}

\section{Subanalytic bundles and sections}

\setcounter{theo}{0}

Let $X\subset M$ a subanalytic set in $M$, with $M$ a real analytic
manifold. Note that unlike analytic subsets, $X$ need not be closed
anymore. We need to define `subanalytic bundles' in some reasonable
fashion, but are hindered by the fact that there is no structure sheaf
for a subanalytic space.

Motivated by the equivalence of (i) and (ii) of Lemma~\ref{analclass} of
the last section, we {\em propose} the following:

\begin{defn}$\left.\right.$\vspace{.5pc}

\noindent {\rm A {\em subanalytic real vector bundle} $E$ {\em of rank}
$k$ {\em generated by an $n$-dimensional vector space} $P$ {\em of
global sections} on a subanalytic set $X\subset M$ is the pullback
$f^{\ast}\nu^{k}$, where $\nu^{k}$ is the universal rank $k$ analytic
quotient bundle on $G(n-k,P)$ defined in Definition \ref{univdef}, and
$f:(X,M)\rarr G(n-k,P)$ is a {\em subanalytic} map (see
Definition~\ref{subanalmapdef}). Then, by pulling back the last two
terms of the exact sequence (\ref{univseq}), there is, by definition, an
exact sequence:
\begin{equation*}
X\times P\stackrel{\eps}{\rarr} E\rarr 0
\end{equation*}
of continuous vector bundles on $X$.\label{subanalbundledef}}
\end{defn}

\begin{rmk}{\rm If $X\subset M$ is an analytic
subset, and $E$ is a real analytic vector bundle of rank $k$, generated
by an $n$-dimensional vector space of real analytic global sections $P$,
we have by (i)$\Rightarrow$ (ii) of Lemma \ref{analclass} that there is
a real analytic classifying map $f:X\rarr G(n-k,P)$. Since $X\subset M$
is an analytic subset, it is subanalytic in $M$, and since the graph
$\Gamma_{f}$ is an analytic space in $M\times G(n-k,P)$, it is
subanalytic, so that $f$ is subanalytic. Thus (ii) of Lemma
\ref{analclass} shows that the definition above is a generalisation of
the real analytic setup to the subanalytic setting.}
\end{rmk}

\begin{exam}{\rm (A typical subanalytic bundle).\ \ Let $X\subset M$ be a
subanalytic set, $M$ a real-analytic manifold, and let $\pi:E\rarr M$ be
a {\em real analytic vector bundle} generated by a vector space $P$ of
global real analytic sections. Thus by definition, there is an exact
sequence of analytic bundles and morphisms:
\begin{equation*}
M\times P\rarr E\rarr 0.
\end{equation*}
By the remark above, there is a real analytic classifying map $f:M\rarr
G(n-k,P)$ such that the exact sequence above is the pullback via $f$ of
the universal sequence
\begin{equation*}
G(n-k,P)\times P\rarr \nu^{k}\rarr 0.
\end{equation*}
If one restricts $f$ to $X$, the restricted map $f_{|X}:(X,M)\rarr G(n-k,P)$
is a subanalytic map, by (v) Proposition~\ref{subanalfacts}. Thus
the restricted bundle $E_{|X}$ (which is the pullback via $f_{|X}$ of
the above sequence on $G(n-k,P)$) is by definition a subanalytic bundle
generated by global sections $P$.\label{typicalexample}}\vspace{-.85pc}
\end{exam}

Again, motivated by the equivalence of (ii) and (iii) of
Lemma~\ref{analclass}, we would like to assert that {\em up to a subanalytic
pseudoequivalence} of $X$ (see Definition~\ref{properequivdef}), the
above Example \ref{typicalexample} of a subanalytic bundle generated by
global sections $P$ is the {\em only} one. More precisely, in complete
analogy with Lemma~\ref{analclass} of the analytic setting in \S3,
we have the following:

\begin{lem} The following are equivalent{\rm :}

\begin{enumerate}
\renewcommand{\labelenumi}{\rm (\roman{enumi})}
\leftskip .1pc
\item $E$ is a rank $k$ subanalytic vector bundle on $X$ generated by an
$n$-dimensional vector space $P$ of global sections in the sense of
Definition~{\rm \ref{subanalbundledef}} above.

\item There exists a real analytic manifold $M_{1}${\rm ,} a subanalytic
map $j:(X,M)\rarr M_{1}$ and a proper map $p:M_{1}\rarr M$ such that{\rm :}

{\rm (a)} $p\circ j= {\rm id}_{X}$, and hence $j:(X,M)\rarr
(j(X),M_{1})$ is a {\em subanalytic pseudoequivalence} in the
sense of Definition~{\rm \ref{properequivdef}}.

{\rm (b)} There is an exact sequence of {\em real analytic
bundles} on the real analytic manifold $M_{1}${\rm :}
\begin{equation*}
\hskip -2.35pc 0\rarr K\rarr M_{1}\times P\stackrel{\eps_{1}}{\rarr} C\rarr 0
\end{equation*}
with $C$ generated by a space of global sections equal to
$P${\rm ,} and such that the pullback of the last two terms via $j$ to $X${\rm ,}
is a surjective morphism of bundles $X\times P\stackrel{\eps}{\rarr} E\rarr 0$.
\end{enumerate}\label{subanallemma}\vspace{-1.5pc}
\end{lem}

\begin{proof}$\left.\right.$

\noindent (i) $\Rightarrow$ (ii).
As in Lemma~\ref{analclass}, set $M_{1}=M\times G(n-k,P)$.
By definition of a subanalytic map, the graph of the subanalytic classifying
map $f$:
\begin{equation*}
X_{1}:=\{(x,y)\in M_{1}: \;y=f(x)\}
\end{equation*}
is a subanalytic space in $M_{1}$. By Example~\ref{graphexample},
the graph embedding:
\begin{align*}
j:(X,M) &\rarr M_{1}\\
x &\mapsto (x,f(x))
\end{align*}
is a subanalytic map. $G(n-k,P)$
is compact, so by Example~\ref{graphexample},
the map $j:(X,M)~\rarr~(X_{1},M_{1})$ is a subanalytic
pseudoequivalence in the sense of Definition~\ref{properequivdef}.

The proof of (b) of (ii) is exactly as in the proof of part (b) of (iii)
in Lemma~\ref{analclass}, and therefore omitted.

\noindent (ii) $\Rightarrow$ (i). From the given exact sequence of real
analytic bundles on $M_{1}$:
\begin{equation}
0\rarr K\longrightarrow M_{1}\times
P\stackrel{\eps_{1}}{\longrightarrow} C\rarr 0 \label{exactm1}
\end{equation}
it follows from (i) and (ii) of Lemma~\ref{analclass} that there is a
real analytic classifying map $g:M_{1}\rarr G(n-k,P)$ (defined by
$g(x)=K_{x}\subset P$) such that the above exact sequence is the
$g$-pullback of the universal exact sequence (\ref{univseq}) on
$G(n-k,P)$. If one composes this classifying map with the the
subanalytic map $j:(X,M)\rarr M_{1}$, which is given to be a subanalytic
pseudoequivalence of $(X,M)$ with $(X_{1},M_{1})$, we have:

\begin{enumerate}
\renewcommand{\labelenumi}{\arabic{enumi}.}
\item $X_{1}=j(X)$ is a subanalytic set in $M_{1}$ by (i) of
Lemma \ref{properlemma}, and

\item The composite $g\circ j:(X,M)\rarr G(n-k,P)$ is a subanalytic map, by
(ii) of Lemma \ref{properlemma} (note that $g:M_{1}\rarr G(n-k,P)$ a
real analytic map implies that $g:(X_{1},M_{1})\rarr G(n-k,P)$ is
subanalytic by (v) of Proposition \ref{subanalfacts}).
\end{enumerate}\vspace{-.5pc}

Clearly this subanalytic composite map:
\begin{equation*}
g\circ j:(X,M)\rarr G(n-k,P)
\end{equation*}
is the classifying map for the bundle $j^{\ast}K$. Thus, letting
$\eps:=j^{\ast}\eps_{1}=j^{\ast}g^{\ast}\phi$, and taking the $j^{\ast}$
of the last two terms of the given exact sequence (\ref{exactm1})
of real analytic bundles on $M_{1}$, we have the exact
sequence
\begin{equation}
X\times P\stackrel{\eps}{\rarr} E\rarr 0 \label{exactx}
\end{equation}
of bundles on $X$, where $E:=j^{\ast}C$. Since the exact sequence
(\ref{exactm1}) is the $g$-pullback of the universal exact sequence
(\ref{univseq}), the above exact sequence (\ref{exactx}) is the pullback
of the last two terms of universal sequence (\ref{univseq}) via the
above subanalytic map $f:=g\circ j:(X,M)\rarr G(n-k,P)$. In particular,
$E$ is a subanalytic bundle of rank $k$ generated by the global sections
$P$, by Definition~\ref{subanalbundledef}. This implies (i).\hfill \ab
\end{proof}

\begin{rmk}
{\rm In complete analogy with Remark~\ref{restrictrmk}, we observe that
(ii) of the above Lemma says that in considering subanalytic bundles on
a subanalytic set $X\subset M$, generated by a vector space $P$ global
sections, we lose no generality (up to subanalytic pseudoequivalence) in
assuming that they are {\em restrictions of real analytic bundles}
(generated by the vector space of analytic sections $=P$) on the {\em
ambient smooth} $M$. \label{restrictrmk2}}
\end{rmk}

\section{Strong Whitney stratifications and transversality}

In this section, we recall some known definitions and results from the
theory of stratifications of subanalytic sets. The general references
for this section are the papers by Verdier [Ve] and Hironaka [Hi].

\setcounter{theo}{0}

\begin{defn}{\rm (The Verdier condition)$\left.\right.$\vspace{.5pc}

\noindent Let $M$ and $M'$ be two locally
closed $C^{\infty}$ submanifolds of some real finite dimensional inner
product space $E$, and such that $M\cap M' =\phi$, with $M'\subset
\ov{M}$.

The {\em property (w)} or {\em Verdier condition} (see \S1.4 in [Ve])
for the pair $(M, M')$ is the following:

For each $y\in \overline{M}\cap M'$, there exists a neighbourhood $U$ of
$y$ in $E$ and a constant $C\in {\Bbb R}^{\ast}_{+}$ such that, for all
$y'\in U\cap M'$ and $x\in U\cap M$, we have
\begin{equation*}
\delta(T_{y'}(M'), T_{x}(M))\leq C \| x -y'\|,
\end{equation*}
where $\delta$ is the distance between two vector subspaces $F$ and $G$
of $E$ (see [Ve], \S1.1) defined by
\begin{equation*}
\delta(F,G) = \sup_{x\in F, \norml x\normr=1}\,\mbox{dist}(x, G),
\end{equation*}
$\mbox{dist}(x, G)$ being the Euclidean distance (in the given norm
$\norml\;\normr$ on $E$) between $x$ and $G$ in $E$ (viz.,
$d(x,G)=\| \pi_{G^{\perp}}(x)\|$). This property is invariant
under smooth local diffeomorphisms of $E$, and hence makes sense for
$M,M'$ contained in a smooth manifold $E$.\label{verdiercond}}
\end{defn}

Now we can define stratifications and {\em Verdier} (or {\em Strong
Whitney} stratifications).

\begin{defn}{\rm (Stratification and Strong Whitney stratification) (see [Ve], \S2.1)$\left.\right.$\vspace{.5pc}

\noindent Let $M$ be a real analytic manifold, countable at $\infty$
(i.e. $M$ is the countable union of compact sets).

A {\em stratification} ${\cal S}$ of $M$ is a partition of $M$ as
$M=\cup_{\alp}M_{\alp}$ satisfying:

\begin{enumerate}
\renewcommand{\labelenumi}{(SW\arabic{enumi})}
\leftskip 1.3pc
\item $M_{\alp}\cap M_{\beta}=\phi$ for $\alp\neq \beta$. Each
$M_{\alp}$ is a locally closed real analytic
submanifold of $M$, smooth, connected, subanalytic in $M$.

\item The family $M_{\alp}$ is locally finite.

\item The family $M_{\alp}$ has the boundary property: i.e.
$\overline{M}_{\alp}\cap M_{\beta}\neq \phi$ implies $M_{\beta}\subset
\overline{M}_{\alp}$.
\end{enumerate}\vspace{-.5pc}

The $M_{\alp}$ are called the {\em strata} of the stratification ${\cal S}$.

A {\em Strong Whitney
stratification}  (or {\em SW-stratification}) ${\cal S}$ of $M$ is a
stratification $M=\cup_{\alp}M_{\alp}$ with the following
additional property:

\begin{enumerate}
\renewcommand{\labelenumi}{(SW\arabic{enumi})}
\setcounter{enumi}{3}
\leftskip 1.3pc
\item If $\overline{M}_{\alp}\supset M_{\beta}$ and if $\alp\neq \beta$,
then the pair $(M_{\alp}, M_{\beta})$ has the property (w) of \ref{verdiercond}
above.
\end{enumerate}

More generally, if $X$ is a subset of $M$, we may define a {\em
stratification} ${\cal S}$ of $X$ to be a partition of $X$ into
$X=\cup_{\alp}X_{\alp}$ where $X_{\alp}$ satisfy (SW1) through (SW3)
above. Similarly, a Strong Whitney or SW-stratification of $X$ to be a
stratification of $X$ having the additional property (SW4). Note that
all the conditions (such as subanalyticity, real analyticity, or the
Verdier condition) {\em are required to hold inside the ambient real
analytic\break manifold} $M$.

Finally, we say a stratification ${\cal S}'$ of $X$ is {\em finer} or a
{\em refinement} of the stratification ${\cal S}$ if each stratum
$X_{\alp}$ of the stratification ${\cal S}$ is a union of strata
$X'_{\beta}$ of the stratification ${\cal S'}$.\label{stratdef}}
\end{defn}

\begin{rmk}$\left.\right.$
{\rm \begin{enumerate}
\renewcommand{\labelenumi}{\rm (\roman{enumi})}
\leftskip .5pc
\item The condition (SW4) above is stronger than Whitney's
condition (b). See Theorem~1.5 (due to Kuo) in [Ve].

\item Since strata are disjoint sets, and constitute a locally finite
collection, a compact set $K\subset M$ can intersect only finitely many
strata. Otherwise, choose a point in $x_{i}\in K\cap M_{i}$ for each
$\alp$ in some countably infinite index set of strata $M_{i}$, and
observe that it will have a limit point $x\in K$ since $K$ is compact.
Every neighbourhood of this $x$ will meet infinitely many $M_{i}$,
contradicting local finiteness. In particular, if $M$ itself is compact,
$M$ is automatically countable at infinity, and the number of strata in
{\em any} stratification of $M$ is finite.

\item If $M$ is SW-stratified as above by strata $\{M_{\alp}\}$,
and $X\subset M$ is any subset which is a union of strata, then $X$ is
also SW-stratified (by the strata of which it is a union).
\end{enumerate}\label{stratrmk}}\vspace{-1pc}
\end{rmk}

The following is a key proposition due to Verdier ([Ve], Theorem~2.2).

\begin{theor}[$($Existence of arbitrarily fine SW-stratifications$)$]
Let $M$ be a real analytic manifold{\rm ,} and $Y_{\beta}$ a locally finite
family of subanalytic subsets of $M$. Then there exists a Strong Whitney
{\rm (}or SW{\rm )}-stratification of $M$ such that each $Y_{\beta}$ is a union of
strata. If $M$ is compact{\rm ,} then the stratification is {\em finite}{\rm ,} i.e.
the number of strata is finite. \label{refine1}
\end{theor}

In particular, for any subanalytic set $X\subset M$, there is a
SW-stratification of $M$ such that $X$ becomes a union of strata, and
hence SW-stratified with the induced stratification. Note that the
assertion for $M$ compact follows from (ii) of Remark~\ref{stratrmk}
above. The analogous statement proving the existence of a Whitney-(b)
stratification for the above setting is Theorem~4.8 of [Hi].

\begin{defn}{\rm (Stratified transversality) (see~1.3.1 on p.~38 in [G-M])$\left.\right.$\vspace{.5pc}

\noindent Let $f:M\rarr N$ be a smooth map of real analytic manifolds.
Let $X\subset M$ and $Y\subset N$ be SW-stratified subsets, with strata
$\{X_{\alp}\}$ and $\{Y_{\bet}\}$ respectively. We say that
$f_{|X}:X\rarr N$ is {\em transverse to} $Y$ (denoted by $f_{|X}\!\transv
Y$) if:
\begin{enumerate}
\renewcommand{\labelenumi}{(\roman{enumi})}
\leftskip .1pc
\item for each stratum $X_{\alp}$ of $X$,
$f_{|X_{\alp}}:X_{\alp}\rarr N$ is a smooth map, and

\item for each stratum $X_{\alp}$ of $X$, and each stratum $Y_{\bet}$
of $Y$, the map $f_{|X_{\alp}}:X_{\alp}\rarr N$ is transverse to
$Y_{\bet}$.
\end{enumerate}\label{transversedef}}\vspace{-1.1pc}
\end{defn}

Note that the above notion does not depend on the analytic, but just the
underlying smooth structures.

\begin{lem} Let $M$ be a real analytic manifold with a SW-stratification
${\cal S}${\rm ,} and let $X\subset M$ be a subanalytic set which is a union
of strata. Let $X_{\alp}$ be the strata of $X${\rm ,} and let ${\cal S}_{X}$
denote this induced stratification of $X$. Let $N\subset M$ be a smooth
real analytic submanifold{\rm ,} subanalytic in $M$. Then if the inclusion map
$i:N\rarr M$ is transverse to the stratification ${\cal S}_{X}${\rm ,} the
intersected stratification ${\cal S}_{X}\cap N$ of $X\cap N$ is defined
by
\begin{equation*}
\{M_{\alp\beta}:\;M_{\alp\beta}\;\mbox{is a connected component
of}\;X_{\alp}\cap N\;\mbox{for some}\;\alp\}.
\end{equation*}

This stratification ${\cal S}_{X}\cap N$ defines a SW-stratification of
$X\cap N$. In particular, each $X_{\alp}\cap N$ is a union of strata
$M_{\alp\bet}$ of ${\cal S}_{X}\cap N$ of the {\em same
dimension}.\label{interstrat}
\end{lem}

\begin{proof}
Since each $X_{\alp}$ is subanalytic and locally closed, and $N$ being a
real analytic submanifold, is also locally closed and given to be
subanalytic, the sets $X_{\alp}\cap N$ are subanalytic and locally
closed, by (i) of Proposition \ref{subanalfacts}. Thus the connected
components $M_{\alp\bet}$ of $X_{\alp}\cap N$ are locally closed, and
subanalytic by (viii) of Proposition~\ref{subanalfacts}. Since
$X_{\alp}$ meets $N$ transversely for each $\alp$, the $M_{\alp\beta}$
are real analytic submanifolds of $M$. Hence (SW1) follows. (SW2) (local
finiteness) also follows from the second statement of (viii) in
\ref{subanalfacts}.

To see (SW3), we need to show that if $M_{\alp\bet}\cap
\ov{M}_{\delta\gamma}\neq \phi$, then $M_{\alp\bet}\subset
\ov{M}_{\delta\gamma}$. Clearly, $\ov{M}_{\delta\gamma}\cap
M_{\alp\bet}$ is closed in $M_{\alp\bet}$. We claim that it is also open
in $M_{\alp\bet}$.

For, let $x\in \ov{M}_{\delta\gamma}\cap M_{\alp\bet}$. This last
intersection being non-empty implies by (SW3) applied to $X_{\alp}$ and
$\ov{X}_{\delta}$, that $x\in X_{\alp} \subset \ov{X}_{\delta}$. By the
topological local triviality of Whitney (b) stratifications (see [G-M],
\S1.4 or [Ma], Corollary~8.4), and the remark in Definition \ref{stratdef}
that SW-stratifications are Whitney (b) stratifications, there exists a
connected neighbourhood $U$ of $x$ in $X_{\alp}$ such that a
neighbourhood $W$ of $x$ in $\ov{X}_{\delta}$ is of the form $W=U\times
N(x)$, where $N(x)$ is the {\em normal slice to $x$} inside
$\ov{X}_{\delta}$. Note that $N(x)$, being the cone on the link $L(x)$,
is connected.

Since the real analytic submanifold $N$ intersects each stratum
transversely, it follows by intersecting everything with $N$, that there
is a neighbourhood $U\cap N$ of $x$ in $X_{\alp}\cap N$ such that the
neighbourhood $W\cap N$ of $x$ in $\ov{X}_{\delta}\cap N$ is of the form
$W\cap N=(U\cap N)\times N(x)$. If we choose $U'$ to be the connected
component of $U\cap N$ containing $x$, and $W'$ to be the connected
component of $W\cap N\cap\ov{M}_{\delta\gam}$ containing $x$, it follows
that $W'$ is of the form $W'=U'\times N'(x)$, where $N'(x)$ is the
intersection of $N(x)$ and $\ov{M}_{\delta\gamma}$. Thus $U'\subset
\ov{M}_{\delta\gamma}\cap M_{\alp\bet}$, and so this last intersection
is open in $M_{\alp\bet}$, and our claim is proven.

The connectedness of $M_{\alp\bet}$ implies that the open and closed
subset $\ov{M}_{\delta\gamma}\cap M_{\alp\bet}$ is equal to
$M_{\alp\bet}$, and (SW3) follows for our stratification ${\cal S}\cap
N$.

To prove the Verdier condition (SW4) for ${\cal S}\cap N$, we first need
the following linear algebraic claim.

\begin{claim} Let $(E, \innp{-}{-})$ be a finite dimensional inner
product space, with a fixed linear subspace $L$. For every subspace $G$
of $E$ intersecting $L$ transversely (viz. $G+L=E$) there exists a
positive constant $C(G)$, depending continuously on $G$, such that for
all linear subspaces $F\subset E$, we have
\begin{equation}
\delta(F\cap L, G\cap L)\leq C(G)\delta(F, G),\label{verdierineq}
\end{equation}
where $\delta$ is the Verdier distance with respect to $\innp{-}{-}$
introduced in Definition~\ref{verdiercond}.
\end{claim}

First note that $L+G=E$ implies $G^{\perp}\cap L^{\perp}=\{0\}$, where
$\bot$ denotes orthogonal complement with respect to $\innp{-}{-}$. Thus
the orthogonal projection $\pi_{G}$ onto $G$ is an isomorphism when
restricted to $L^{\perp}$. Hence the subspace $\pi_{G}(L^{\perp})$ has
dimension $=\dim\,L^{\perp}=\mbox{codim}\,L$ which, by transversality of
$G$ and $L$, is precisely the codimension of $G\cap L$ in $G$. Also, for
$y\in G\cap L,\;z\in L^{\perp}$, we have
$\innp{y}{\pi_{G}(z)}=\innp{\pi_{G}(y)}{z}=\innp{y}{z}=0$, thus implying
that $\pi_{G}(L^{\perp})$ is orthogonal to $G\cap L$. Hence
$L_{1}:=\pi_{G}(L^{\perp})$ is the orthogonal complement of $G\cap L$ in
$G$.

Define a new inner product $\innp{-}{-}'$ on $E$ by setting
$\innp{y}{z}'=\innp{y}{z}$ for $y,z\in L$, $\innp{y}{z}'=\innp{y}{z}$
for $y,z\in L_{1}$, and $\innp{L_{1}}{L}'=0$. The orthogonal complement
of a subspace $W$ with respect to $\innp{-}{-}'$ will be denoted
$W^{\top}$. Then, with respect to this new inner product $\innp{-}{-}'$,
we have an orthogonal direct sum decomposition
\begin{equation*}
E=L_{1}\oplus' L,
\end{equation*}
where `$\oplus '$' signifies that the decomposition is orthogonal with
respect to $\innp{-}{-}'$. Also, $L^{\top}=L_{1}\subset G$. Thus
$G^{\top}\subset L$. Hence, by counting dimensions, and noting that
$\innp{G\cap L}{G^{\top}}'=0$, it follows that
\begin{equation*}
L= (G\cap L)\oplus' G^{\top}.
\end{equation*}

Hence for any $x\in F\cap L$, we will have
$\pi_{G^{\top}}(x)=\pi_{(G\cap L)^{\top}}(x)$. Thus, denoting the
Verdier distance with respect to $\innp{-}{-}'$ by $\delta'$, we have
\begin{align}
\delta'(F\cap L, G\cap L) &= \sup_{\norml x\normr'=1, x\in F\cap L}
\| \pi_{(G\cap L)^{\top}}(x)\|'\nonumber\\
&= \sup_{\norml x\normr'=1, x\in F\cap L}\| \pi_{G^{\top}}(x)\|'\nonumber\\
&\leq \sup_{\norml x\normr'=1, x\in F}\| \pi_{G^{\top}}(x)\|'\nonumber\\
&= \delta'(F,G).\label{deltaprime}
\end{align}
(All orthogonal projections in the last six lines are with respect to
$\langle -,-\rangle'$.)

Now, by the above definition of $\innp{-}{-}'$, it follows that there is
a positive constant $A(G)$ depending continuously on $G$ such that for
every $x\in E$:
\begin{equation*}
\fr{1}{A(G)}\norml x\normr \leq \norml x\normr ' \leq A(G)\norml x \normr
\end{equation*}
from which it is easy to deduce that there is an inequality
\begin{equation*}
\delta'(F,G)\leq C(G)\delta(F,G),
\end{equation*}
where $C(G)$ is a positive constant depending continuously on $G$. Also
since $\innp{-}{-}$ and $\innp{-}{-}'$ agree on $L$, it follows that
$\delta'(F\cap L, G\cap L)=\delta(F\cap L, G\cap L)$. Hence the
inequality (\ref{verdierineq}) follows by applying (\ref{deltaprime}),
and hence the Claim.

If $y\in M_{\alp\bet}\subset \ov{M}_{\delta\gam}$, then $y\in X_{\alp}$
and we have by (SW4) applied to $y\in X_{\alp}\subset \ov{X}_{\delta}$,
that there is some neighbourhood $U$ of $y$ such that for all $y'\in
U\cap X_{\alp}$ and $x\in X_{\delta}\cap U$ we have
\begin{equation*}
\delta(T_{y'}(X_{\alp}), T_{x}(X_{\delta}))\leq C\| x -y'\|.
\end{equation*}
Let $U'$ be the connected component of $U\cap M_{\alp\bet}$ containing
$y$, and simultaneously locally trivialise on $U$ the manifold pair
$(M,N)$ so that the bundle pair $(TM,TN)_{|U}$ is isomorphic to $U\times
(E,L)$. Set $F:=T_{y'}(X_{\alp})$ and $G:=T_{x}(X_{\delta})$, so that
$T_{y'}(M_{\alp\bet})=F\cap L$ and $T_{x}(M_{\delta\gamma})=G\cap L$, by
the fact that $N$ meets $X_{\alp}$ and $X_{\delta}$ transversely. For
$x\in U$, one can find a bound $C$ for the constant $C(G)$ in
(\ref{verdierineq}) by the continuity of $C(G)$ in $G$. Now we apply the
inequality (\ref{verdierineq}) above to get the Verdier condition (w) on
$U'$. Thus (SW4) is verified for the stratification ${\cal S}\cap N$.
This proves the lemma.\hfill \ab
\end{proof}

\begin{rmk}{\rm We also need a straighforward extension of
Proposition~\ref{interstrat}. Namely, in the above situation let $N$ be
a closed connected subset, which is a real analytic submanifold with
real analytic boundary $\boun N$. This can be regarded as a
SW-stratified subset of $M$ with two strata, viz. $N^{\circ}:=N\setminus
\boun N$ and $\boun N$. For this stratification, (SW1), (SW2) and
(SW3) are straightforward, and (SW4) holds because at any point $x\in
\boun N$, the germ of the triple $(N,\boun N, M)$ is (by definition)
real analytically isomorphic to $(\R^{k-1}\times\R_{+}, \R^{k-1},
\R^{n})$, which obviously satisfies (SW4). Let $X\subset M$ be a
SW-stratified subset, such that $X_{\alp}\!\transv N^{\circ}$ and
$X_{\alp}\!\transv \boun N$ for each stratum $X_{\alp}$ of $X$. Then the
intersection stratification on $X\cap N$ defined by taking the connected
components of all the intersections $X_{\alp}\cap N^{\circ}$ and
$X_{\alp}\cap \boun N$ is an SW-stratification of $X\cap N$. The proof
is exactly the same as in the case of \ref{interstrat} above.
\label{interstrat2}}
\end{rmk}

\begin{defn}{\rm (Stratified submersion) (see [G-M], \S1.5 on p.~41 and [Ve], 3.2)$\left.\right.$\vspace{.5pc}

\noindent Let $f:M\rarr N$ be a smooth map, and let $X\subset M$ be a
SW-stratified subset, with strata $\{X_{\alp}\}$. We say $f$ is a {\em
stratified submersion} if for each stratum $X_{\alp}$ of $X$, the
restriction $f_{|X_{\alp}}:X_{\alp}\rarr N$ is smooth and a submersion.
\label{stratsubdef}}
\end{defn}

Again, the above notion depends only on the underlying smooth
structures. In [Ve], Definition~3.2, such a map is referred to as a map
{\em `transverse to the stratification'} (on $X$).
The following is a key proposition due to Verdier.

\begin{theor}[$($First Isotopy Lemma$)$] (See ${[\hbox{\it Ve}]}$, Theorem~{\rm 4.14}).\ \
Suppose that $X$ is a {\em closed} SW-stratified subset in $M${\rm ,} with
stratification ${\cal S}$. Let $N$ be a smooth real analytic manifold{\rm ,}
$f:M\rarr N$ be a real analytic map{\rm ,} {\em proper on} $X${\rm ,} and a
stratified submersion. Let $y_{0}\in N${\rm ,} and $M_{0}$ and $X_{0}$ be the
fibres over $y_{0}$ of $M$ and $X$ respectively. {\rm (}Note that by the
hypothesis{\rm ,} $M_{0}$ meets all the strata of $X$ transversely{\rm ,} so that
$X_{0}=X\cap M_{0}$ acquires an induced stratification ${\cal S}_{0}$
via the connected components of $X\cap M_{0}${\rm ,} as in {\rm \ref{interstrat}}
above{\rm )}. Then there exists an open neighbourhood $V$ of $y_{0}$ in $N${\rm ,}
and a homeomorphism of the stratified spaces $(X\cap f^{-1}(V), {\cal
S})$ onto $(X_{0}\times V, {\cal S}_{0}\times V)$ preserving the
stratifications and compatible with the projections to $V$. {\rm (}Again note
that $X\cap f^{-1}(V)${\rm ,} being an open subset of $X${\rm ,} has a natural
induced stratification{\rm ,} also denoted by ${\cal S}${\rm ,} from $X${\rm ,} and
$X_{0}\times V$ has a natural product stratification ${\cal S}_{0}\times
V$.{\rm )}\label{firstiso}
\end{theor}

In fact, in Theorem~4.14 of [Ve], Verdier proves that the homeomorphism
above is a `rugeux' (coarse) homeomorphism, which is slightly stronger
than saying that it is a homeomorphism. We do not need this stronger
statement. The corresponding statement for Whitney (b) stratifications
is Thom's First Isotopy Lemma, and due to Thom and Mather (see [7]
in [Ve]).

\section{A tubular neighbourhood theorem for subanalytic bundles}

Let $X$ be a topological space, and $E$ be a continuous real vector
bundle of rank $k$ on $X$. Let $\norml\,\normr$ be some continuous
bundle metric, and for $\eps>0$ let $E(\eps)$ denote the $\eps$-disc
bundle of $E$ with respect to this bundle metric. Denote the
zero-section of this bundle by $0_{E}$.

\setcounter{theo}{0}

\begin{defn}{\rm (Tubular neighbourhoods)$\left.\right.$\vspace{.5pc}

\noindent Let $\pi:E\rarr X$ be as above and let $s:X\rarr E$ be a
continuous section. We will say that $s$ {\em has a tubular
neighbourhood in} $X$ if there exists a neighbourhood $V$ of
$s^{-1}(0_{E})$ in $X$, and $\eps >0$, and a homeomorphism
\begin{equation*}
\Phi :V\rarr E(\eps)_{|s^{-1}(0_{E})}
\end{equation*}
such that the composite
\begin{equation*}
s^{-1}(0_{E})\hookrightarrow V\stackrel{\Phi}{\rarr}
E(\eps)_{|s^{-1}(0_{E})}
\end{equation*}
is the map $x\rarr s(x)=0_{x}$ for all $x\in s^{-1}(0_{E})$.
\label{tubnbddef}}
\end{defn}

\begin{rmk}$\left.\right.$

{\rm \begin{enumerate}
\renewcommand{\labelenumi}{(\roman{enumi})}
\leftskip .1pc
\item If $E\rarr X$ and $F\rarr X$ are isomorphic real vector bundles,
via a continuous bundle isomorphism $\tau:E\rarr F$, then clearly a
section $s$ of $E$ will have a tubular neighbourhood in $X$ iff the
section $\tau\circ s$ of $F$ has a tubular neighbourhood in $X$. We just
need to observe that $s^{-1}(0_{E})=(\tau\circ s)^{-1}(0_{F})$, and that
$\tau$ will induce a homeomorphism of the disc bundle $F(\eps)$ with the
disc bundle $E(\eps)$ where $E$ is given the $\tau$-pullback bundle
metric from $F$.

\item More generally, if there is a bundle diagram:
\begin{align*}
\begin{array}{lll}
E_{1} &\stackrel{\tau}{\longrightarrow} &E_{2}\\
\raise .2pc\hbox{$_{\pi_{1}}$} \downarrow& &\downarrow \raise .2pc\hbox{$_{\pi_{2}}$}\\
X_{1}&\stackrel{f}{\longrightarrow} &X_{2}
\end{array}
\end{align*}
with $f$ a {\em homeomorphism}, and $\tau$ a continuous vector bundle
equivalence, then the section $s_{2}$ of $E_{2}$ has a tubular
neighbourhood in $X_{2}$ iff the section $f^{\ast}s_{2}:=\tau^{-1}\circ
s_{2}\circ f$ of $E_{1}$ has a tubular neighbourhood in $X_{1}$. This is
clear because we will have a continuous bundle equivalence $E_{1}\rarr
f^{\ast}E_{2}$ on $X$, the statement (i) above, and the fact that the
homeomorphism $f$ induces a homeomorphism of any neighbourhood $V$ of
$s_{2}^{-1}(0_{E_{2}})$ with the neighbourhood $f^{-1}(V)$ of
$s_{1}^{-1}(0_{E_{1}})$.
\end{enumerate}\label{tubnbdrmk1}}\vspace{-.6pc}
\end{rmk}

Before proving the main result, we need some lemmas.

\begin{lem}{\rm (}Transversality to the $0$-section in smooth case{\rm )} {\rm (}see
Theorem~{\rm 1.3.6} on p.~{\rm 39} in {\rm [G-M])}.\ \ Let $M$ be a real
analytic manifold{\rm ,} and $\pi:E\rarr M$ be a real analytic vector bundle
of rank $k$ on $M${\rm ,} generated by an $n$-dimensional vector space $P$ of
real analytic global sections. Let $Y\subset M$ be a real analytic
submanifold {\rm (}in particular{\rm ,} locally closed in $M${\rm )}. Then there exists a
Baire subset {\rm (}i.e. whose complement is of measure $0${\rm ,} in particular{\rm ,}
dense{\rm )} $U\subset P$ such that $s\in U$ implies that $s_{|Y}:Y\rarr E${\rm ,}
is transverse to the zero-section $0_{E}$ of $E$.\label{transmooth}
\end{lem}

\begin{proof}
By definition, and (ii) of Lemma~\ref{analclass}, we have the exact
sequence of real analytic bundles on $M$:
\begin{equation*}
M\times P\stackrel{\eps}{\rarr} E\rarr 0
\end{equation*}
which is the pullback via the analytic classifying map $f:M\rarr
G(n-k,P)$ of the last two terms of the universal exact sequence
(\ref{univseq}) on $G(n-k,P)$. In particular $\eps=f^{\ast}\phi$ above
is real analytic. The section $s\in P$ is viewed as an analytic section
of $E$ by taking the map $\eps_{s}=\eps(-,s):M\rarr E$ of $E$. Since
$\eps$ is a smooth epimorphism of bundles, it is a smooth submersion.

Restricting the above real analytic map $\eps$ to $Y$, we have a real
analytic map:
\begin{equation*}
Y\times P\stackrel{\theta}{\rarr} E.
\end{equation*}
We claim that this map $\theta$ is transverse to the zero-section
$0_{E}$ of $E$.

\pagebreak

For, if $(y,s)\in Y\times P$ such that $\theta(y,s)=0_{y}\in E_{y}$, the
fact that $\theta_{y}=\theta(y,-)$ is a linear surjection shows that the
partial derivative $\partial_{s}\theta_{y}=\theta_{y}$ takes the linear
subspace $0\oplus T_{s}P=0\oplus P$ of $T_{(y,s)}(Y\times P)$
surjectively onto the tangent space $T_{0_{y}}(E_{y})=E_{y}$ of the
fibre $E_{y}$. Hence the image of the total derivative $D\theta(y,s)$
contains $T_{0_{y}}(E_{y})=E_{y}$. Since
$T_{0_{y}}(E)=T_{0_{y}}(0_{E})\oplus E_{y}$, it follows that
$T_{0_{y}}(0_{E})$ and $\im{D\theta(y,s)}$ together span $T_{0_{y}}(E)$.
Hence $\theta$ is transverse to the zero-section $0_{E}$ of $E$.

Since $0_{E}$ is a codimension $k$ smooth submanifold of $E$, the
inverse image $W:=\theta^{-1}(0_{E})$ is a smooth real-analytic
submanifold of $Y\times P$, of dimension $\dim\,W=\dim\,Y +n -k$. Let
$p:Y\times P\rarr P$ denote the second projection, and consider the
following\break diagram:
\begin{align*}
\begin{array}{rccc}
W\hookrightarrow &Y\times P&\stackrel{\theta}{\rarr} &E_{|Y}\\
p_{|W}\searrow &\downarrow\,{_{p}} & &\\
 & & &\\
 &P&  &
\end{array}.
\end{align*}

Note that the fibre $W_{s}=p^{-1}(s)\cap W=p_{|W}^{-1}(s)$ is precisely
the set of zeroes of the section $\theta_{s}=\eps_{s|Y}=s_{|Y}$ of the restricted bundle
$E_{|Y}\rarr Y$.  We make
the following:

\begin{claim}
The section $\theta_{s}=\theta(-,s)=s:Y\rarr E$ is transverse to the
zero-section $0_{E}$ of $E$ iff $s\in P$ is a regular value of
$p_{|W}:W\rarr P$.
\end{claim}

\begin{poclaim} Let
\begin{equation*}
q:T_{0_{y}}(E)\rarr T_{0_{y}}(E)/T_{0_{y}}(0_{E})=E_{y}
\end{equation*}
denote the natural quotient map.

Since $\theta:Y\times P\rarr E$ is transverse to
$0_{E}$, for each $(y,s)\in W$ we have the diagram:
\begin{align*}
\begin{array}{ll@{\ \ }l@{\ \ }lll@{\ }l}
&\ker{Dp_{|W}(y,s)}&\longrightarrow &T_{y}(Y)&
\stackrel{q\circ D\theta_{s}(y)}{\longrightarrow}&
T_{0_{y}}(E)/T_{0_{y}}(0_{E})&\\
&\downarrow & &\downarrow &&\downarrow &\\
0\rarr &T_{(y,s)}(W) &\longrightarrow &T_{(y,s)}(Y\times P)
&\stackrel{q\circ D\theta(y,s)}{\longrightarrow}
&T_{0_{y}}(E)/T_{0_{y}}(0_{E})\rarr &0\\
&\downarrow {_{Dp_{|W(y,s)}}} & &\downarrow {_{Dp}} & &\downarrow &\\
0\rarr &P &\longrightarrow &P &\stackrel{\qquad\quad\ \ }\longrightarrow  &0\rarr 0
\end{array}
\end{align*}
in which the columns and last two rows are exact. By the Snake Lemma, and
the fact that $\coker{Dp}=0$, we
have an exact sequence:
\begin{equation*}
0\rarr \ker{Dp_{|W}}\rarr T_{y}(Y)
\stackrel{q\circ D\theta_{s}(y)}{\longrightarrow}
E_{y}\rarr\coker{Dp_{|W}}\rarr 0.
\end{equation*}

Consequently, the map $q\circ D\theta_{s}(y)$ is surjective iff
$Dp_{|W}(y,s)$ is surjective. Since
$p_{|W}^{-1}(s)=\theta_{s}^{-1}(0_{E})~\times~\{s\}$, $s$ is a regular
value of $p_{|W}$ iff $q\circ D\theta_{s}(y)$ is surjective for all
$y\in \theta_{s}^{-1}(0_{E})$, that is, iff $\theta_{s}:Y\rarr E$ is
transverse to $0_{E}$. Hence the Claim.

By Sard's theorem, there is a Baire subset $U\subset P$ of regular
values for $p_{|W}:W\rarr P$. Hence the lemma follows, for this choice
of $U$.\hfill \ab
\end{poclaim}
\end{proof}

We have the following immediate corollary for subanalytic sets.

\begin{coro}{\rm (}Stratified transversality to the $0$-section{\rm )}$\left.\right.$\vspace{.5pc}

\noindent Let $M$ be a real analytic manifold which is countable at
$\infty$. Let ${\cal S}$ be a SW-stratification of $M${\rm ,} and let
$\{M_{\alp}\}_{\alp\in\Lambda}$ denote the strata {\rm (}see
Definition~{\rm \ref{stratdef}}{\rm )}. Let $E$ be a real analytic vector bundle of rank $k$ on
$M${\rm ,} generated by an $n$-dimensional real vector space $P$ of
real-analytic global sections. Then there exists a Baire {\rm (}in particular
dense{\rm )} subset $U\subset P$ such that for $s\in U${\rm ,} we have
\begin{enumerate}
\renewcommand{\labelenumi}{{\rm (\roman{enumi})}}
\leftskip .5pc
\item $s:M\rarr E$ is transverse to $0_{E}${\rm ,} and

\item $s_{|M_{\alp}}\!\transv 0_{E}$ {\rm (}see Definition~{\rm \ref{transversedef})}
for each $\alp\in \Lambda$.

\item For any $X\subset M$ which is a union of strata {\rm (}and hence
itself SW-stratified{\rm ),} $s_{|X}\!\transv 0_{E}$.
\end{enumerate}\label{transubanal}\vspace{-1pc}
\end{coro}

\begin{proof} Since the manifold $M$ is a countable union of compact
sets, and a compact set can intersect only finitely many members of
locally finite collection, and SW-stratifications are locally finite, it
follows that the SW-stratification above consists of only
countably many strata. We thus assume that $\Lambda={\Bbb N}$.
In particular, $M$ is the union:
\begin{equation*}
M=\cup_{i=1}^{\infty}M_{i}
\end{equation*}
of countably many strata. Set $M_{0}=M$, as a notational convention.

By (i) of Definition~\ref{stratdef}, each $M_{i}$ is a locally closed
real analytic submanifold of $M$. Applying the previous Lemma
\ref{transmooth} to $Y=M_{i}$ for $i=0, 1, ...$, we find a Baire subset
$U_{i}\subset P$ such that for $s\in U_{i}$, $s_{|M_{i}}\!\transv
0_{E}$. The set $U=\cap_{i=0}^{\infty}U_{i}$, being the countable
intersection of Baire sets, is a Baire set (the countable union of
measure zero sets is measure zero). Clearly for $s\in U$, (i) and (ii)
follow. (iii) is obvious from (ii) by\break
Definition~\ref{transversedef}.\hfill \ab
\end{proof}

\begin{theor}[$($Stratified Tubular Neighbourhood Theorem~1$)$] Let $M$
be a {\em compact} real analytic manifold{\rm ,} and $X\subset M$ a
subanalytic set {\rm (}not necessarily closed{\rm )} in $M$. In accordance with
Theorem~{\rm \ref{refine1}} above{\rm ,} equip $M$ with a SW-stratification ${\cal
S}$ by strata $\{M_{i}\}_{i=1}^{m}$ such that the subanalytic set $X$ is
a union of strata. {\rm (}Note that the strata form a finite collection by
{\rm (ii)} of Remark~{\rm \ref{stratrmk})}. Let $\pi:E\rarr M$ be a real
analytic vector bundle of rank $k${\rm ,} generated by an $n$-dimensional
vector space $P$ of real analytic global sections. Then{\rm ,} there exists a
Baire subset $U\subset P$ such that for $s\in P${\rm ,} the section $s_{|X}$
of the restricted subanalytic bundle $E_{|X}\rarr X$ has a tubular
neighbourhood in the sense of\break Definition~{\rm \ref{tubnbddef}}.
\label{strattubular}
\end{theor}

\begin{proof} The proof is rather involved, though the essential idea is
contained in the sketch of the proof of Theorem~1.11, p. 47 in [G-M].

The main steps are as follows: Equip the compact manifold $M$ with a
finite (SW)-stratification by $M_{i}$ so that the subanalytic set $X$ is
a union of strata, by Theorem \ref{refine1}. By Corollary
\ref{transubanal}, for a Baire subset $U\subset P$ and $s\in U$, the
restriction of $s$ to {\em each stratum} is transverse to $0_{E}$. Let
$Z$ denote the zero locus of $s$ in $M$, and $Z_{i}=Z\cap M_{i}$. Next
we choose a thin enough $\eps$-disc bundle in $E$ which intersects only
those strata $s(M_{i})$ of $s(M)$ that meet $0_{E}$ (i.e. for which
$Z_{i}\neq\phi$). Now, shrink down this $\eps$-disc bundle to the
zero-section $0_{E}$ via the scaling map $e\mapsto te$ of $E$, $t\in
\R$. The final step is to prove, using the First Isotopy Lemma
\ref{firstiso} that the $t$-family of intersections of these shrinking
neighbourhoods with $s(M)$ is a stratified product $Y\times I$, where
$I$ is some interval containing $0$. For this family, the fibre over
$0\in I$ is the bundle restricted to $Z$, and the fibre over $a\neq 0\in
I$, is a homeomorphic image of the above neighbourhood. Thus these
fibres are both stratified homeomorphic to $Y$. Because it is a
stratified product, and $X$ is a union of strata, one can restrict the
(stratification preserving) homeomorphisms above to $X$. We carry out
the details below.

Define $M_{0}=M$ for notational convenience. We know by
Corollary \ref{transubanal} that for the given stratification of
$M=\cup_{i=1}^{m}M_{i}$ above, there exists a Baire subset $U\subset P$
such that $s\in U$ implies that $s_{i}:=s_{|M_{i}}$ is transverse to the
zero section $0_{E}$ of $E$ for each $i=0,1,...,m$. (Note that
$s_{i}:M_{i}\rarr E_{|M_{i}}$ is transverse to the zero section of
$E_{|M_{i}}$ iff $s_{i}:M_{i}\rarr E$ is transverse to the zero section
$0_{E}$ of $E$.) In particular $s_{|X}$ is transverse to the zero
section $0_{E}$ of $E$ (or, in the notation of Definition
\ref{transversedef}, \hbox{$s_{|X}\!\transv 0_{E}$}, for the real analytic map
$s:M\rarr E$). Thus the zero-locus $Z:=Z_{0}:=s^{-1}(0_{E})$ is a real
analytic submanifold of $M$, and the restriction of $Z$ to $M_{i}$, viz.
the zero-locus $Z_{i}=Z\cap M_{i}=s_{i}^{-1}(0_{E})$ is a real analytic
submanifold of $M_{i}$, and hence $M$, for each $i=1,...,m$. Denote the
restricted bundle $E_{|M_{i}}$ by $E_{i}$.

We shall prove that for such an $s\in U$, the section $s_{|X}:X\rarr
E_{|X}$ has a tubular neighbourhood in $X$, in the sense of Definition
\ref{tubnbddef}.

Since $\pi:E\rarr M$ is isomorphic as a real analytic bundle to
$f^{\ast}\nu^{k}$, and since $\nu^{k}\rarr G(n-k,P)$ is a real analytic
subbundle of $G(n-k,P)\times P$ via the real analytic splitting coming
from a constant (hence real analytic) bundle metric on $G(n-k,P)\times
P$, it follows that $\pi:E\rarr M$ has a real analytic pullback bundle
metric $\norml\;\normr$. Let $p:M\times \R\rarr M$ denote the first
projection. Then the bundle $\pi_{1}:=(\pi\times 1):E\times \R\rarr
M\times \R$, which is precisely $p^{\ast}E$, has the pulled back metric,
also denoted $\norml\;\normr$, also real analytic. Note that $M\times
\R$ has the product stratification denoted ${\cal S}\times \R$, and
$X\times \R$ is a union\break of strata.

Similarly $E$ has the SW-stratification ${\cal S}_{E}=\pi^{-1}({\cal
S})$, with strata $E_{i}=E_{|M_{i}}= \pi^{-1}(M_{i})$. (It is easy to
verify, using the local analytic triviality of $E$, that $\pi^{-1}({\cal
S})$ satisfies (SW1) through (SW4). Also $E_{i}$ are subanalytic inside
$E$ because $M_{i}$ are subanalytic in $M$ and (vi) of Proposition
\ref{subanalfacts}.) Similarly, $E\times \R$ is SW-stratified by ${\cal
S}_{E}\times \R$, which is the same as the SW-stratification
$(\pi_{1})^{-1}({\cal S}\times\R)$.

Consider the real analytic map:
\begin{align*}
\psi:E\times\R &\rarr E\\
(e,t) &\mapsto te
\end{align*}
and denote $\psi(-,t):E\rarr E$ by $\psi_{t}$. Denote the restricted
bundle map $\psi_{|E_{i}\times \R}\rarr E_{i}$ by $\psi_{i}$.

\def\claaimm{\trivlist\item[\hskip\labelsep{\it Claim {$(1)_{0}$}.}]}

\begin{claaimm}
The map $\psi:E\times \R\rarr E$ is transverse to the smooth
real analytic submanifold $s(M)\subset E$.
\end{claaimm}

In fact we shall show that for each $t\in \R$, the map $\psi_{t}:E\rarr
E$ is transverse to $s(M)$. Since for each $t\in\R$, the image
$\im{D\psi}(e,t)$ contains $\im{D\psi_{t}(e)}$, this will prove the
assertion. If $t\neq 0$, we have $\psi_{t}:E\rarr E$ is a real analytic
diffeomorphism, and hence transverse to $s(M)$. If $t=0$, we have
$\psi_{0}:E\rarr E$ equal to projection onto the zero-section $0_{E}$,
viz. the map $e\mapsto 0_{\pi(e)}$. Thus for $e\in E$, if we denote
$\pi(e)=x\in M$, we have $\psi_{0}(e)=0_{x}$ and
$\im{D\psi_{0}(e)}=T_{0_{x}}(0_{E})$. If $\psi_{0}(e)\in s(M)$, it
follows that $e=s(x)=0_{x}$, so that $x\in s^{-1}(0_{E})$. By the choice
of $s$, $s$ is transverse to $0_{E}$, so that we have
$\im{Ds(x)}+T_{0_{x}}(0_{E})=T_{0_{x}}(E)$. That is,
$T_{s(x)}(s(M))+\im{D\psi_{0}(e)}= T_{0_{x}}(E)$. Thus $\psi_{0}$ is
also transverse to $s(M)$, and our claim follows.

By the transversality above, it follows that $N:=\psi^{-1}(s(M))$ is a
real analytic submani- fold of $E\times\R$, whose fibre over $t\in\R$ is
the real analytic manifold $N(t):=\psi_{t}^{-1}(s(M))$. Since $s:M\rarr
E$ is a real analytic embedding, and sections of bundles are proper
maps, it follows that $s(M)$ is a smooth real analytic subspace of $E$.
In particular it is closed. Since $\psi$ is real analytic, it follows
that $N=\psi^{-1}(s(M))$ is a real analytic subspace in $E\times \R$
which is closed. By the transversality of $\psi$ to $s(M)$, it follows
that $N$ is a smooth real analytic subspace of $E\times \R$.

Next we have

\def\claaimm{\trivlist\item[\hskip\labelsep{\it Claim {$(2)_{0}$}.}]}

\begin{claaimm}
The real analytic submanifold $N$ is transverse to the zero section $0_{E\times\R}$, and the intersection
$N\cap\,0_{E\times\R}$ is $s(Z)\times\R$, where $Z:=s^{-1}(0_{E})$ is the
zero locus of $s$.
\end{claaimm}

Since $s:M\rarr E$ is transverse to $0_{E}$ by the choice of
$s$ above, it follows that for each $0_{x}\in s(M)$ ($\Leftrightarrow
x\in Z=s^{-1}(0_{E})$) we have
\begin{equation}
T_{0_{x}}(s(M))+T_{0_{x}}(0_{E})=T_{0_{x}}(E).\label{tangeqn1}
\end{equation}

Now for $(0_{x},t)\in N\cap\,0_{E\times\R}$, we have the identifications
$T_{0_{x}}(E)=T_{0_{x}}(0_{E})+E_{x}$ and $T_{(0_{x},t)}(E\times
\R)=T_{0_{x}}(0_{E})+E_{x}+\R$. With this identification, clearly
$D\psi(0_{x},t)(W,0,0)=W$ for $W\in T_{0_{x}}(0_{E})$, since $\psi$
restricted to $0_{E}\times \R$ is the first projection to $0_{E}$. Let
us denote $L:=D\psi(0_{x},t)$ for notational convenience. Then by the
foregoing remark and (\ref{tangeqn1}) above, we have:
\begin{equation}
T_{0_{x}}(s(M)) + L(T_{(0_{x},t)}(0_{E\times\R}))=T_{0_{x}}(E).
\label{tangeqn2}
\end{equation}

For any linear map $L$, any two linear subspaces $A,B$, we
have the identity $L^{-1}A+B=L^{-1}(A+L(B))$. Substituting $A=
T_{0_{x}}(s(M))$ and $B=T_{(0_{x},t)}(0_{E\times\R})$, and using
(\ref{tangeqn2}) above we find
\begin{align*}
L^{-1}T_{0_{x}}(s(M)) + T_{(0_{x},t)}(0_{E\times\R}) &=
L^{-1} (T_{0_{x}}(s(M)))+ L(T_{(0_{x},t)}(0_{E\times\R}))\\
&= L^{-1}(T_{0_{x}}(E))\\
&= T_{(0_{x},t)}(E\times\R).
\end{align*}

Since by Claim $(1)_{0}$,  $\psi$ is transverse to $s(M)$,
we have $L^{-1}T_{0_{x}}(s(M))=T_{(0_{x},t)}(N)$. Thus we have
\begin{equation*}
T_{(0_{x},t)}(N)+T_{(0_{x},t)}(0_{E\times\R})=
T_{(0_{x},t)}(E\times\R).
\end{equation*}
That is, $N$ meets $0_{E\times \R}$ transversely. The first part of
Claim $(2)_{0}$ follows.

To see the second part of Claim $(2)_{0}$, note that  $(0_{x},t)\in
N\cap 0_{E\times \R}$ iff $\psi(0_{x},t)=0_{x}$ is in $s(M)$. But
$0_{x}\in s(M)$ iff $0_{x}=s(x)$ iff $x\in s^{-1}(0_{E})=Z$.
Thus $N\cap 0_{E\times \R}=s(Z)\times \R$, and Claim $(2)_{0}$ is proved.

In the proofs of the above two claims, the transversality of $s$ to the
zero section was the only fact used. By the choice of $s\in U$ we have
$s_{i}:M_{i}\rarr E_{i}$ is transverse to the zero section $0_{E_{i}}$
for all $i$. Thus, in the proofs of the above two claims, we replace $M$
by the stratum $M_{i}$, $E$ by $E_{i}$, $s$ by $s_{i}:=s_{|M_{i}}$,
$\psi$ by $\psi_{i}:= \psi_{|E_{i}\times\R}$, and repeat everything
(noting that $s_{i}:M_{i}\rarr E_{i}$ is transverse to $0_{E_{i}}$ iff
$s_{i}:M_{i}\rarr E$ is transverse to $0_{E}$). By so doing, for
$i=0,1,...,m$, we obtain:

\def\claaim{\trivlist\item[\hskip\labelsep{\it Claim {$(1)_{i}$}.}]}

\begin{claaim}
$\psi_{i}:E_{i}\times \R\rarr E_{i}$ is transverse to $s(M_{i})$, and thus
$N_{i}:=\psi_{i}^{-1}(s(M_{i}))=\psi^{-1}(s(M_{i}))$ is a smooth
real analytic subspace of $E_{i}$, and is equal to the intersection
$N\cap (E_{i}\times \R)$. (For notational consistency, we set $N_{0}:=N$.)
\end{claaim}

\def\claaim{\trivlist\item[\hskip\labelsep{\it Claim {$(2)_{i}$}.}]}

\begin{claaim}
$N_{i}:=\psi^{-1}(s(M_{i}))$ intersects $0_{E\times\R}$ transversely, and
this intersection is $s(Z_{i})\times \R=s(Z\cap M_{i})\times\R$.
\end{claaim}

From the given SW-stratification ${\cal S}$ in the hypothesis of
the theorem, the smooth analytic subspace $s(M)$ gets an
induced SW-stratification ${\cal S}_{s}$, defined by the strata
$s(M_{i}),\; i=1,2,...,m$. (Note $s(M_{i})$ are subanalytic in $E$ by
(ii) of Proposition \ref{subanalfacts} because $s$, being a
continuous section of a bundle, is proper.)

\def\claaim{\trivlist\item[\hskip\labelsep{\it Claim {$(3)$}.}]}

\begin{claaim}
The manifold $s(M)$ meets the strata $E_{i}$ of ${\cal S}_{E}$
transversely (inside $E$), and $N$ meets the strata $E_{i}\times\R$ of
${\cal S}_{E\times\R}$ transversely (inside $E\times\R$).
\end{claaim}

To prove the first assertion, let $s(x)\in s(M)\cap E_{i}=s(M_{i})$. Since
$\pi:E\rarr M$ satisfies $\pi\circ s =\mbox{id}_{M}$, the
image $\im{Ds(x)}=T_{s(x)}(s(M))$ is a vector space complement to
$\ker{D\pi(s(x))}=T_{s(x)}(E_{x})=E_{x}$ inside $T_{x}(E)$, and hence
\begin{equation}
T_{s(x)}(s(M))+ E_{x}=T_{s(x)}(E).\label{smtrans}
\end{equation}

Since $E_{x}\subset T_{s(x)}(E_{i})$ for all $x\in M_{i}$, it
follows that
\begin{equation*}
T_{s(x)}(s(M))+ T_{s(x)}(E_{i})=T_{s(x)}(E)
\end{equation*}
and hence $s(M)$ is transverse to $E_{i}$ for all $i$. The first
assertion of Claim (3) follows.

Now we prove the second assertion of Claim (3), i.e. $N$ meets
$E_{i}\times\R$ transversally for all $i=1,2,...,m$ inside $E\times \R$.
We have already observed that $D\psi(e,t):T_{(e,t)}(E\times\R)\rarr
T_{te}(E)$ is surjective for all $t\neq 0$ (because $\psi_{t}:E\rarr E$
is a diffeomorphism). Also for such a $t\neq 0$, $D\psi(e,t)$ maps
$E_{x}\oplus\R$ onto $E_{x}$. Hence, from (\ref{smtrans}) we conclude
that
\begin{equation}
T_{s(x)}(s(M)) + D\psi(e,t)(E_{x}\oplus\R)=T_{s(x)}(E)\label{smtrans0}
\end{equation}
for $(e,t)\in N\cap (E_{i}\times\R)$ with $t\neq 0$, and $\pi(e)=x$.
Applying the identity $L^{-1}(A+L(B))=L^{-1}A+B$ to the linear map
$L=D\psi(e,t),\;A=T_{s(x)}(s(M)),\; B=E_{x}\oplus\R$, and noting that
$L^{-1}(T_{s(x)}(s(M))=T_{(e,t)}(N)$ by definition and Claim~$(1)_{0}$,
we obtain
\begin{equation*}
T_{(e,t)}(N)+(E_{x}\oplus\R)=T_{(e,t)}(E\times\R).
\end{equation*}
Since $E_{x}\oplus\R=T_{(e,t)}(E_{x}\times\R)\subset
T_{(e,t)}(E_{i}\times \R)$ for all $i=0,1,...,m$, we obtain that the
intersection of $N$ with $E_{i}\times \R$ is transverse at all points
$(e,t)$ for $t\neq 0$.

Now we look at the situation where $t=0$. If a point $(e,0)\in N\cap
(E_{i}\times\R)$, we have that $\psi(e,0)=0_{x}\in s(M_{i})=s(M)\cap
E_{i}$ and $s(x)=0_{x}$. From the paragraph preceding the proof of Claim
$(1)_{i}$, we have that $s:M_{i}\rarr E_{i}$ is transverse to the
zero-section $0_{E_{i}}$ in $E_{i}$, so that
\begin{equation}
T_{0_{x}}(s(M_{i}))+T_{0_{x}}(0_{E_{i}})=T_{0_{x}}(E_{i}).\label{smtrans2}
\end{equation}

Since $Ds(x)(T_{x}(M_{i}))$ is a complement to $E_{x}$ in $T_{x}(E_{i})$, for each $x\in M_{i}$, and
$i=1,2,...,m$, it follows that
\begin{equation}
T_{0_{x}}(E_{i})=E_{x}+Ds(x)(T_{x}(M_{i})).\label{smtrans3}
\end{equation}

Let $\nu_{i}$ denote the normal bundle to $M_{i}$ in $M$, and
$\nu_{i,x}$ its fibre at $x\in M_{i}$. Since
$T_{0_{x}}(s(M))=Ds(x)(T_{x}(M))$ and
$T_{0_{x}}(s(M_{i}))=Ds(x)(T_{x}(M_{i}))$, we have
\begin{align}
T_{s(x)}(s(M))\!+\! T_{0_{x}}(0_{E_{i}}) &= Ds(x)(T_{x}(M))+
T_{0_{x}}(0_{E_{i}})\nonumber\\
&= Ds(x)(T_{x}(M_{i}))+Ds(x)(\nu_{i,x})+T_{0_{x}}(0_{E_{i}})\nonumber \\
&= T_{0_{x}}(E_{i}) + Ds(x)(\nu_{i,x})\;\;\;\;\;\;\;\;(\mbox{by eq. (\ref{smtrans2})})\nonumber\\
&= E_{x}\!+\!Ds(x)(T_{x}(M_{i}))\!+\!Ds(x)(\nu_{i,x})\ \ \ (\mbox{by eq. (\ref{smtrans3})})\nonumber\\
&= E_{x}+Ds(x)(T_{x}(M))\nonumber\\
&= T_{0_{x}}(E).\label{smtrans4}
\end{align}

$\psi_{0}:E\rarr E$ is just the projection of $E$ to its zero
section $0_{E}$, and hence
$D\psi_{0}(e)(T_{(e,0)}(E_{i}))$ $=T_{0_{x}}(0_{E_{i}})$, so that
$D\psi(e,0)(T_{(e,0)}(E_{i}\times\R))$ contains $T_{0_{x}}(0_{E_{i}})$.
Combining with eq.~(\ref{smtrans4}), we have
\begin{equation*}
T_{s(x)}(s(M))+D\psi(e,0)(T_{(e,0)}(E_{i}\times\R))=T_{0_{x}}(E).
\end{equation*}

Now we apply our lemma $L^{-1}A + B=L^{-1}(A+L(B))$ to the linear map
$L=D\psi(e,0)$, $A= T_{s(x)}(s(M))$, $B = T_{(e,0)}(E_{i}\times\R)$,
and noting that by definition and Claim $(1)_{0}$,
$L^{-1}(T_{s(x)}(s(M))=T_{(e,0)}(N)$, we obtain
\begin{equation*}
T_{(e,0)}(N)+ T_{(e,0)}(E_{i}\times\R)=T_{(e,0)}(E\times\R)
\end{equation*}
which shows that $N$ meets $E_{i}\times\R$ transversely at $(e,0)$
and our Claim (3) follows.

Since $s(M_{i})=E_{i}\cap s(M)$, the first assertion of Claim (3) shows that
the stratification ${\cal S}_{s}$ is just the intersection stratification
${\cal S}_{E}\cap s(M)$ as defined in Lemma~\ref{interstrat}. In particular,
by that same lemma, ${\cal S}_{s}$ is a (SW)-stratification of $s(M)$.

By the second assertion of Claim (3) and Lemma \ref{interstrat}, we can
give the SW-stratification ${\cal S}_{N}:=
{\cal S}_{E\times\R}\cap N$ to the smooth
real analytic subspace $N=\psi^{-1}(s(M))$. The strata are precisely the
connected components of $N_{i}=N\cap (E_{i}\times\R)=\psi^{-1}(s(M_{i}))$.

Now we need another assertion:

\def\claaim{\trivlist\item[\hskip\labelsep{\it Claim {$(4)$}.}]}

\begin{claaim}
For each $i=0,1,...,m$, define $p_{i}:N_{i}\rarr \R$ to be
the restriction of the second projection $p:E\times\R\rarr \R$ to
$N_{i}$.  The derivative
\begin{equation*}
Dp_{i}(0_{x},t):T_{(0_{x},t)}(N_{i})\rarr \R
\end{equation*}
is surjective at all points $(0_{x},t)\in N_{i}\cap
0_{E\times\R}=s(Z_{i})\times\R$.
\end{claaim}

We now prove Claim (4). For the curve $s\mapsto (0_{x},t+s)$ in
$E\times\R$, starting at $(0_{x},t)$ with initial velocity $(0,0,1)\in
T_{(0_{x},t)}(E\times\R)=T_{0_{x}}(0_{E})+E_{x}+\R$, we have
$\psi(0_{x},s+t)=(s+t)0_{x}\equiv 0_{x}$ for all $s$ so that
\begin{equation*}
D\psi(0_{x},t)(0,0,1)= 0
\end{equation*}
for all $t\in \R$. Thus the subspace $0\oplus 0\oplus \R\subset
T_{(0_{x},t)}(E\times\R)$ lies in $\ker{D\psi(0_{x},t)}$. Since (by
Claim $(1)_{i}$ above), $T_{(0_{x},t)}(N_{i})=
(D\psi(0_{x},t))^{-1}(T_{0_{x}}s(M_{i}))$, it contains this kernel
$\ker{D\psi(0_{x},t)}$. Thus
\begin{equation*}
T_{(0_{x},t)}(N_{i})\supset 0\oplus 0\oplus\R
\end{equation*}
for all $(0_{x},t)\in N_{i}$. Since $Dp_{i}$ maps the subspace $(0\oplus
0\oplus\R)$ isomorphically onto $\R=T_{t}(\R)$, it follows that
\begin{equation*}
Dp_{i}(0_{x},t):T_{(0_{x},t)}(N_{i})\rarr \R
\end{equation*}
is surjective at all
points $(0_{x},t)\in N_{i}\cap 0_{E\times\R}$ and all $i=0,1,...,m$.
Hence Claim (4) is proved.

Recall that $Z_{i}:=Z\cap M_{i}$, and by the second part of Claim
$(2)_{i}$, we have $N_{i}\cap 0_{E\times\R}=s(Z_{i})\times\R$.

From Claim $(4)$, and by continuity of $Dp_{i}$, we have the following claim.

\def\claaim{\trivlist\item[\hskip\labelsep{\it Claim {$(5)$}.}]}

\begin{claaim}
For each $i$ such that $Z_{i}\neq\phi$, and for $(0_{x},t)\in N_{i}\cap
0_{E\times\R}=s(Z_{i})\times\R$, there exists a neighbourhood $U_{x,t}$
of $(0_{x},t)$ in $E\times\R$ such that the derivative
\begin{equation*}
Dp_{i}(e,s):T_{(e,s)}(N_{i})\rarr\R
\end{equation*}
is surjective for all $(e,s)\in U_{x,t}\cap N_{i}$.
\end{claaim}

Renumber the strata so that $Z_{i}:=Z\cap M_{i}\neq \phi$ for
$i=1,2,...,r$ and $Z_{i}=\phi$ for $i=r+1,...,m$. Again, for notational
convenience, set $Z_{0}:=Z=s^{-1}(0_{E})$. For $\delta >0$, let
$E(\delta)$ and $S(\delta)$ denote the open $\delta$-disc bundle and the
$\delta$-sphere bundles of $E$ respectively (with respect to the above
real analytic bundle metric $\norml\;\normr$).

\def\claaim{\trivlist\item[\hskip\labelsep{\it Claim {$(6)$}.}]}

\begin{claaim}
For each $i=0,1,...,r$, and for $(0_{x},t)\in
N_{i}\cap\,0_{E\times\R}=s(Z_{i})\times\R$, there exists a neighbourhood
$U_{x,t}$ of $(0_{x},t)$ in $E\times\R$ such that $N_{i}\cap U_{x,t}$
meets $(S(\delta)\times \R)\cap U_{x,t}$ transversely for all $\delta
>0$. (Note that for $\delta$ very large, the intersection
$(S(\delta)\times \R)\cap U_{x,t}$ might be empty, in which case this
assertion is vacuously true.)
\end{claaim}

To see this, we first prove that for $i=0,1,...,r$ and
$0_{x}\in s(M_{i})\cap 0_{E}=s(Z_{i})$, there is a neighbourhood $U_{x}$
of $0_{x}$ in $E$ such that $s(M_{i})\cap U_{x}$ meets $S(\delta)\cap U_{x}$
transversely for all $\delta>0$. Since this assertion is local around
$0_{x}$, and $E$ is locally trivial, we may assume without loss of
generality that $E$ is trivial.

By this triviality, there is a natural linear surjection:
\begin{equation*}
\rho(e):T_{e}(E)\rarr E_{\pi(e)}=(\pi^{\ast}E)_{e}
\end{equation*}
for each $e\in E$. The fact that $s(M_{i})$ meets $0_{E}$ transversely
at $0_{x}$ implies that the composite
\begin{equation*}
T_{0_{x}}(s(M_{i}))\hookrightarrow
T_{0_{x}}(E)\stackrel{\rho(0_{x})}{\rarr} (\pi^{\ast}E)_{0_{x}}=E_{x}
\end{equation*}
is surjective. By continuity, there exists a neighbourhood $U_{x}$ of
$0_{x}$ in $E$ such that the composite
\begin{equation*}
T_{e}(s(M_{i}))\hookrightarrow
T_{e}(E)\stackrel{\rho(e)}{\rarr} (\pi^{\ast}E)_{e}
\end{equation*}
is surjective for all $e\in s(M_{i})\cap U_{x}$.

At a point $e\in S(\delta)$, it is obvious that the one-dimensional
quotient space $T_{e}(E)/T_{e}(S(\delta))$ is a quotient of
$E_{\pi(e)}=(\pi^{\ast}E)_{e}$ (by the tangent space to the sphere fibre
$S(\delta)_{\pi(e)}$). Thus the composite map
\begin{equation*}
T_{e}(s(M_{i}))\hookrightarrow
T_{e}(E)\stackrel{\rho(e)}{\rarr} (\pi^{\ast}E)_{e}\rarr
T_{e}(E)/T_{e}(S(\delta))
\end{equation*}
is also surjective for all $e\in s(M_{i})\cap S(\delta)\cap U_{x}$. But
this is precisely the statement that $s(M_{i})\cap U_{x}$ meets
$S(\delta)\cap U_{x}$ transversely inside $E$.

For the sake of convention, we define $S(\delta)=0_{E}$ for $\delta =
0$. Then, by the fact that $s(M_{i})$ is transverse to $0_{E}$, we have
that $s(M_{i})\cap U_{x}$ meets $S(\delta)\cap U_{x}$ transversely in
$E$ for all $\delta\geq 0$.

Now let $t\in\R$, and let $U_{x,t}$ be a neighbourhood of $(0_{x},t)$
such that $\psi(U_{x,t})\subset U_{x}$, $U_{x}$ as above. Hence
$(v,\lam)\in (S(\delta)\times\R)\cap N_{i}\cap U_{x,t}$ implies that
$e:=\psi(v,\lam)=\lam v\in (S(|\lam|\delta))\cap s(M_{i})\cap U_{x}$.

By the above choice of the neighbourhood $U_{x}$, at the point
$e=\lam v\in U_{x},\;\lam\in\R$, we have
\begin{equation}
T_{e}(S(|\lam|\delta))+T_{e}(s(M_{i}))=T_{e}(E).\label{tangineq3}
\end{equation}

If $(v,\lam)\in U_{x,t}$ with $\lam\neq 0$, then $\psi_{\lam}:E\rarr E$
is a diffeomorphism. In this event, the derivative
\begin{equation*}
D\psi(v,\lam):T_{(v,\lam)}(E\times\R)\rarr T_{\lam v}(E)
\end{equation*}
is surjective. Taking the inverse image under this surjective map
$D\psi(v,\lam)$ of both sides of the equality (\ref{tangineq3}), noting
that $\psi$ is transverse to $s(M_{i})$, and that
$\psi_{\lam}:S(\delta)\times \{\lam\}\rarr S(|\lam|\delta)$ is a
diffeomorphism for $\lam\neq 0$, we have
\begin{equation*}
T_{(v,\lam)}(S(\delta)\times \{\lam\}) + T_{(v,\lam)}(N_{i})=
T_{(v,\lam)}(E\times\R).
\end{equation*}
Since $T_{(v,\lam)}(S(\delta)\times \R)$ contains the subspace
$T_{(v,\lam)}(S(\delta)\times \{\lam\})$, we see that at all points
$(v,\lam)$ with $\lam\neq 0$, which lie in
$U_{x,t}$ and in $(S(\delta)\times\R)\cap N_{i}$,
the intersection of $S(\delta)\times
\R$ and $N_{i}$ is transverse.

To check the case $\lam =0$, let $(v,0)\in (S(\delta)\times\R)\cap
N_{i}\cap U_{x,t}$. Note $\psi(v,0)=0_{y}$ where $y=\pi(v)$ and
$\pi:E\rarr M$ is the bundle projection. Since $(v,0)\in N_{i}\cap
U_{x,t}$, we have $0_{y}\in s(M_{i})\cap U_{x}$. Thus, setting
$L:=D\psi(v,0)$, we have $L(T_{(v,0)}(S(\delta)\times
\{0\}))=T_{0_{y}}(0_{E})$, since the projection $\pi:S(\delta)\rarr M$ is a
submersion. Since $0_{y}\in s(M_{i})\cap 0_{E}$, and $s(M_{i})$
intersects $0_{E}$ transversely, we have
\begin{equation*}
T_{0_{y}}(s(M_{i}))+T_{0_{y}}(0_{E})=T_{0_{y}}(E)
\end{equation*}
from which it follows that
\begin{equation*}
T_{0_{y}}(s(M_{i}))+ L(T_{(v,0)}(S(\delta)\times
\{0\}))=T_{0_{y}}(E)
\end{equation*}
which implies
\begin{equation*}
T_{0_{y}}(s(M_{i}))+ L(T_{(v,0)}(S(\delta)\times
\R))=T_{0_{y}}(E).
\end{equation*}
Again we use the identity $L^{-1}(A+L(B))=L^{-1}A+B$ for any linear map
$L$, and any subspaces $A,\,B$ and obtain, by setting
$A=T_{0_{y}}(s(M_{i}))$ and $B=T_{(v,0)}(S(\delta)\times\R)$, that
\begin{equation*}
L^{-1}T_{0_{y}}(s(M_{i}))+ T_{(v,0)}(S(\delta)\times
\R)=T_{(v,0)}(E\times \R)
\end{equation*}
which is to say,
\begin{equation*}
T_{(v,0)}(N_{i})+T_{(v,0)}(S(\delta)\times
\R)=T_{(v,0)}(E\times \R).
\end{equation*}
Thus the intersection of $N_{i}\cap U_{x,t}$ and $(S(\delta)\times\R)\cap
U_{x,t}$ is transverse. This proves Claim~$(6)$.

\pagebreak

Since $M$ is compact, $0_{E}$ is a compact subset of $E$. Since, by
Claim $(2)_{i}$ above, $N_{i}$ intersects $0_{E\times\R}$ transversely for
each $i=0,1,...,m$, each connected component of each $N_{i}$ intersects
$0_{E\times\R}$ transversely inside $E\times\R$. Thus no such connected
component can be contained in $0_{E\times\R}$. That is, no stratum of
the stratification ${\cal S}_{N}$ of the smooth manifold
$N=\psi^{-1}(s(M))$ is contained in $0_{E\times\R}$.

We now note that since every stratum $s(M_{i})$ of $s(M)$ meets $0_{E}$
transversely, $s(M_{i})~\subset~\ov{s(M_{j})}=s(\ov{M}_{j})$ for
$i\neq j$ will imply that the codimension of $s(M_{i})\cap 0_{E}$ in
$\ov{s(M_{j})}\cap 0_{E}$ is the same as the codimension of $s(M_{i})$
in $\ov{s(M_{j})}$, which is non-zero. Thus if $\ov{s(M_{j})}\cap 0_{E}\neq \phi$, the
union
\begin{equation*}
\bigcup\left\{s(M_{i})\cap 0_{E}:\;s(M_{i})\subset \ov{s(M_{j})},\;\;
\mbox{and}\;i\neq j\right\}
\end{equation*}
has non-zero codimension in $\ov{s(M_{j})}\cap 0_{E}$. Thus
\begin{equation*}
\hskip -4pc s(M_{j})\cap 0_{E}=\left(\ov{s(M_{j})}\cap 0_{E}\right)
\Big\backslash \left(\bigcup\left\{s(M_{i})\cap 0_{E}:\;s(M_{i})\subset
\ov{s(M_{j})}\;\;\mbox{and}\;i\neq j\right\}\right)\neq \phi.
\end{equation*}$\left.\right.$\vspace{-1.4pc}

That is, $\ov{s(M_{j})}\cap 0_{E}\neq \phi$ if and only if
$s(M_{j})\cap 0_{E}\neq \phi$, if and only if $j=0,1,...,r$.

By noting that $N_{j}=\psi^{-1}(s(M_{j}))$, and
$\psi(0_{E\times\R})=0_{E}$, the same fact obtains for the $N_{j}$, viz.
$\ov{N}_{j}~\cap~0_{E\times\R}\neq \phi$ if and only if $N_{j}\cap
0_{E\times\R}\neq \phi$ which, in turn, happens if and only if
$j=0,1,...,r$.

Consider the closed subset $C$ of $E\times\R$ defined by
\begin{equation*}
C:=\cup_{i=r+1}^{m} \ov{N}_{i}
\end{equation*}
which is disjoint from $0_{E\times\R}$, by the preceding paragraph.
Since $C$ is closed and $0_{E}\times [-2,2]$ is compact, there will
exist an $\eps > 0$ such that the restricted $2\eps$-disc bundle
$E(2\eps)\times [-2,2]$ does not intersect $C$. Thus $E(2\eps)\times [-
2,2]$ is disjoint from $N_{i}$ for $i=r+1,...,m$. We saw above in
Claim $(2)_{i}$ that $N_{i}\cap 0_{E\times\R}=s(Z_{i})\times\R$, and similarly
$N_{i}\cap (0_{E}\times [-2,2])=s(Z_{i})\times [-2,2]$. Thus $N_{i}$
intersects $E(2\eps)\times [-2,2]$ if and only if $N_{i}$ intersects
$0_{E}\times [-2,2]$, and this happens if and only if $Z_{i}\neq \phi$,
i.e. if and only if $i=0,1,...,r$. We record this fact in claim (7)

\def\claaim{\trivlist\item[\hskip\labelsep{\it Claim {$(7)$}.}]}

\begin{claaim}
$(E(2\eps)\times [-2,2])\cap N_{i}=\phi$ iff $i=r+1,r+2,...,m$.
\end{claaim}

By reducing $\eps$ if necessary, and from Claims (5) and (6),
using the compactness of $0_{E}~\times~[-2,2]$, we can further
assert the following Claim~(8).

\def\claaim{\trivlist\item[\hskip\labelsep{\it Claim {$(8)$}.}]}

\begin{claaim}
For $i=0,1,...,r$, the derivative $Dp_{i}(e,s)$ is surjective
for all  $(e,s)\in \left(E(2\eps)\times [-2,2]\right)\cap N_{i}$.
\end{claaim}

\def\claaim{\trivlist\item[\hskip\labelsep{\it Claim {$(9)$}.}]}

\begin{claaim}
For $i=0,1,...,r$, the intersection of
$S(\delta)\times [-2, 2]$ with $N_{i}$ is transverse for all
$\delta <2\eps$.
\end{claaim}

We need an analogue of Claim (8) above for $S(\eps)$. That is

\def\claaim{\trivlist\item[\hskip\labelsep{\it Claim {$(10)$}.}]}

\begin{claaim}
For $\eps$ small enough, and $0\leq i\leq r$ such
that the intersection\break \hbox{$(S(\eps)\times[-2,2])\cap
N_{i}\neq \phi$,} the derivative $Dp_{i}(e,s)$ is surjective
for all  $(e,s)\in (S(\delta)\times [-2,2])\cap N_{i}$ and all $\delta
<2\eps$.
\end{claaim}

We note that the tangent subspace $T_{(e,s)}(S(\delta)\times \R)$ has
the line complement $\R e\subset E_{\pi(e)}$ in the tangent space
$T_{(e,s)}(E\times \R)=T_{(e,s)}(E(2\eps)\times\R)$. Since $Dp$
annihilates all vectors in $T_{e}(E(2\eps))\oplus \{0\}\subset
T_{(e,s)}(E)$, it annihilates the subspace $E_{\pi(e)}$, and hence $\R
e$. Thus the $Dp_{i}$ image of $T_{(e,s)}((S(\delta)\times [-2,2])\cap
N_{i})$ is the same as the $Dp_{i}$ image of $T_{(e,s)}((E(2\eps)\times
[-2,2])\cap N_{i})$, which is all of $\R=T_{s}([-2,2])$ by Claim (5)
above. Thus Claim (10) is proved.

Consider $\ov{E(\eps)}\times (-2,2)$ as a SW-stratified subspace of the
analytic manifold $E\times (-2,2)$ with just the two strata
$E(\eps)\times (-2,2)$ and $S(\eps)\times (-2,2)$, as in Remark
\ref{interstrat2}. By the fact that $E(\eps)\times (-2,2)$ is open
(and subanalytic in $E \times (-2,2)$, since it is the inverse image of
$(-\eps^{2},\eps^{2})$ under the analytic map $\norml\;\normr^{2}$), and
the fact that $N$ is SW-stratified by the connected components of
$N_{i}$, it follows that $A:=N\cap (E(\eps)\times (-2,2))$ is
SW-stratified by the connected components of the analytic submanifolds
$N_{i}\cap (E(\eps)\times (-2,2))$. By Claim (7) above, we have
\begin{equation*}
A=\cup_{i=1}^{r}(N_{i}\cap (E(\eps)\times (-2,2)).
\end{equation*}

Similarly by Claim (9) above, and Remark~\ref{interstrat2}, the subset
\hbox{$\boun A=N\cap (S(\eps)\times (-2,2))$} is SW-stratified by the
connected components of the analytic submanifolds\break \hbox{$N_{i}\cap
(S(\eps)\times (-2,2))$.} (Note that $S(\eps)$ is a smooth
real analytic subspace of $E$ since $\norml\;\normr^{2}$ is a real
analytic function, and $\ov{E}(\eps)$ is a real analytic
manifold with boundary $S(\eps)$). In fact, by Remark~\ref{interstrat2},
\begin{equation*}
\ov{A}:=N\cap (\ov{E}(\eps)\times (-2,2))
\end{equation*}
is SW-stratified by the connected components of the real analytic
submanifolds $N_{i}\cap (S(\eps)\times (-2,2))$ and $N_{i}\cap
(E(\eps)\times (-2,2))$ for $i=1,...,r$. Let us call these connected
components $A_{\alp}$, where $\alp\in F$ for some finite set $F$. Thus
the subset $\ov{A}$ is a SW-stratified space in $E\times (-2,2)$ with
stratification ${\cal S}_{A}$ by $A_{\alp}$.

By Claims (8) and (10) above, for $i=0,1,...,r$ the maps
$p_{i}:N_{i}\cap (S(\eps)\times (-2,2))\rarr (-2,2)$ and
$p_{i}:N_{i}\cap (E(\eps)\times (-2,2))\rarr (-2,2)$ are submersions.
Thus they are submersions when restricted to each connected component.
In particular, $p_{i}:A_{\alp}\rarr (-2,2)$ is a submersion for each
$\alp\in F$. That is $p:\ov{A}\rarr (-2,2)$ is a stratified submersion
in the sense of Definition \ref{stratsubdef}.

For any compact subset $K\subset (-2,2)$,
$p_{|\ov{A}}^{-1}(K)=p^{-1}(K)\cap \ov{A}$ is a closed set, and
contained in the compact set $\ov{E}(\eps)\times K$. ($\ov{E}(\eps)$ is
compact since $M$ is compact!). Thus $p_{|\ov{A}}^{-1}(K)$ is compact,
implying that $p_{|\ov{A}}$ is proper.

By the First Isotopy Lemma \ref{firstiso}, applied to the analytic map
$p:E\times (-2,2)\rarr (-2,2)$ and the closed SW-stratified subset
$\ov{A}\subset E\times (-2,2)$, it follows that for each point $t\in
(-2,2)$, there is a neighbourhood $U_{t}:=(t-\delta, t+\delta)$ of $t$,
and a stratum preserving rugeux homeomorphism
\begin{equation*}
h:\ov{A}\cap p^{-1}(U_{t})\rarr \ov{A}_{t}\times (t-\delta,t+\delta),
\end{equation*}
where $\ov{A}_{t}:=\ov{A}\cap (E\times \{t\})$, such that
$pr_{2}\circ h=p$ (here $pr_{2}$ is the second projection on the right
hand side). That is, $\ov{A}\rarr (-2,2)$ is a topologically
locally-trivial stratified fibre bundle.

By the compactness and connectedness of $[0,1]\subset (-2,2)$, there is
therefore a stratification preserving rugeux homeomorphism
\begin{equation*}
h: \ov{A}_{1}\rarr\ov{A}_{0}
\end{equation*}
between the fibres $\ov{A}_{1}$ and $\ov{A}_{0}$. But $\ov{A}_{1}= N\cap
(\ov{E}(\eps)\times\{1\})=(\psi_{1}^{-1}(s(M))\cap
(\ov{E}(\eps))\times\{1\})$ by definition. Since $\psi_{1}:E\rarr E$ is
the identity map, this last set is stratified homeomorphic to $(s(M)\cap
\ov{E}(\eps))$. On the other hand, the map $\psi_{0}:E\rarr E$ is the
bundle projection map $e\rarr 0_{\pi(e)}$. Thus
$\psi_{0}^{-1}(s(M))=E_{|Z}$ where $Z=s^{-1}(0_{E})$ is the zero-locus
of $s$. Thus $\ov{A}_{0}$ is stratified homeomorphic to
$\ov{E}(\eps)_{|Z}$.

Note that $s(M)\cap E(\eps)$ is an open neighbourhood of $s(M)\cap
0_{E}=s(Z)\cap 0_{E}$. By the last para, this neighbourhood is
stratified homeomorphic to $E(\eps)_{|Z}$. The fact that this
homeomorphism preserves strata shows that it maps $s(M_{i})\cap E(\eps)$
homeomorphically to $E(\eps)_{|Z_{i}}$, where $Z_{i}=Z\cap M_{i}$ for
$i=1,2,...,r$ (recall for $i>r$, the above sets are empty). Since $X$ is
a union of strata from among the $M_{i}$, and this last homeomorphism is
stratum preserving, it follows that this homeomorphism when restricted
to the open neighbourhood $s(X)\cap E(\eps)_{|X}$ of $s(X)\cap
0_{E}=s_{|X}(Z\cap X)$, is again a stratum preserving homeomorphism
$s(X)\cap E(\eps)_{|X}\rarr E(\eps)_{|Z\cap X}$.

Since the analytic embedding $s:M\rarr s(M)$ defines a homeomorphism
between a neighbourhood of $s^{-1}(0_{E})\cap X$ and $s(X)\cap
E(\eps)_{|X}$, Theorem \ref{strattubular} follows. \hfill \ab
\end{proof}

\begin{theor}[$($Stratified Tubular Neighbourhood Theorem~2$)$]
Let $M$ be a compact real analytic manifold{\rm ,} and $X\subset M$ be a
subanalytic set. Let $E$ be a subanalytic bundle of rank $k$ over $X$
generated by an $n$-dimensional vector space of global sections $P$ in
the sense of Definition {\rm \ref{subanalbundledef}}. Then{\rm ,} for a Baire subset
$U\subset P${\rm ,} and $s\in U,$ there exists a neighbourhood of $Z:=s^{-
1}(0_{E})$ in $X$ which is homeomorphic to an $\eps$-disc bundle
$E(\eps)_{|Z}$ for some $\eps >0$ {\rm (}i.e. $s:X\rarr E_{|X}$ has a tubular
neighbourhood in the sense of Definition~{\rm \ref{tubnbddef})}.
\label{strattubular2}
\end{theor}

\begin{proof}
By Lemma \ref{subanallemma}, there is a subanalytic
pseudoequivalence $j:(X,M)\rarr (X_{1},M_{1})$, with $X_{1}$ a
subanalytic set in $M_{1}$, such that our given bundle $E$ is the
pullback $j^{\ast}C$ of a real analytic rank-$k$ vector bundle $C$ on
$M_{1}$, generated by a space of global analytic sections $P$. Further,
since $M$ is compact, and $M_{1}=M\times G(n-k,P)$, $M_{1}$ is also
compact.

By (ii) of Remark~\ref{tubnbdrmk1}, the Tubular Neighbourhood
Theorem for $E$ on $X$ follows from the Stratified Tubular
Neighbourhood Theorem~1 (i.e. Theorem~\ref{strattubular}) applied to
$C$, $X_{1}$ and $M_{1}$ above. \hfill \ab
\end{proof}

\begin{rmk}
{\rm Theorem~\ref{strattubular2} covers the case of $X$ being any real
projective or affine variety. In fact, any real projective or affine
algebraic constructible set can be regarded as a subanalytic subset in
projective space. In the analytic situation too, any real analytically
constructible subset of a compact real analytic manifold automatically
becomes a subanalytic set. The main Theorem \ref{strattubular2} applies
in all of the above situations, provided the bundle $E$ is a subanalytic
bundle generated by global sections in the sense of Definition
\ref{subanalbundledef}. \label{strattubular3}}
\end{rmk}

\section*{Acknowledgements}

I would like to thank V~Srinivas for asking me the question addressed in
this paper, and many fruitful discussions. I also thank Jishnu Biswas
for carefully reading this paper and pointing out errors, and for many
fruitful conversations on subanalytic sets. Finally, I thank the referee
for valuable comments and corrections.

\end{document}